\documentclass[letterpaper,conference, 10pt]{ieeeconf}
\IEEEoverridecommandlockouts
\usepackage{amsmath,graphicx,epsfig,color,amsfonts,subfigure,algorithm, psfrag,parskip}
\usepackage{version,xspace}
\usepackage[noend]{algorithmic}
\graphicspath{{./fig/}}

\newcommand{\algostep}[1]{{\small\texttt{#1:}}\xspace}

\def\prob{\mathbb{P}}

\def\expt{\mathbb{E}}
\def\real{\mathbb{R}}

\def\naturals{\mathbb{N}}
\def\ones{\boldsymbol{1}}
\def\t{\boldsymbol{t}}
\def\m{\boldsymbol{\mu}}
\def\trans{\text{T}}

\def\e{\boldsymbol{e}}
\renewcommand{\natural}{\mathbb{N}}

\newcommand{\until}[1]{\{1,\dots, #1\}}

\newcommand{\subscr}[2]{#1_{\textup{#2}}}

\newcommand{\setdef}[2]{\{#1 \; | \; #2\}}

\newcommand{\map}[3]{#1: #2 \rightarrow #3}
\newcommand{\union}{\operatorname{\cup}}

\newcommand{\subject}{\text{subject to}}
\newcommand{\maximize}{\text{maximize}}
\newcommand{\minimize}{\text{minimize}}

\newcommand{\argmax}{\operatorname{argmax}}

\newcommand\oprocendsymbol{\hbox{$\square$}}
\newcommand\oprocend{\relax\ifmmode\else\unskip\hfill\fi\oprocendsymbol}

\renewcommand{\baselinestretch}{1}

\renewcommand{\theenumi}{\roman{enumi}}

\newtheorem{theorem}{Theorem}
\newtheorem{lemma}[theorem]{Lemma}
\newtheorem{remark}[theorem]{Remark}

\newtheorem{example}[theorem]{Example}

\newcommand{\Tproc}{\subscr{\mathcal{T}}{processed}}
\def\tinf{t_{\text{inf}}}

\def\fdi{f'_{\text{inv}}}
\title{Task Release Control for Decision Making Queues
\thanks{This work has been supported in part by AFOSR MURI Award
    FA9550-07-1-0528.} }
\author{Vaibhav~Srivastava \hspace{0.5in} Ruggero~Carli \hspace{0.5in}Francesco~Bullo\hspace{0.5in} C{\'e}dric Langbort \\ 
\thanks{Vaibhav Srivastava  and Francesco Bullo are with the Center for Control, Dynamical Systems, and   Computation, University of California, Santa Barbara, Santa Barbara, CA
93106, USA,\tt{ \{vaibhav,bullo\}@engineering.ucsb.edu}} 
\thanks{Ruggero Carli is with the Department of Information Engineering, University of Padova ,  Padova, Italy, \tt{carlirug@dei.unipd.it
}}
\thanks{C{\'e}dric Langbort is with the Coordinated Science Laboratory, 
University of Illinois at Urbana-Champaign, Urbana, IL 61801, USA, \tt{langbort@illinois.edu}}}
\begin{document}
\maketitle
\begin{abstract}
We consider the optimal duration allocation in a decision making queue. Decision making tasks arrive at a given rate to a human operator. The correctness of the decision made  by human evolves as a sigmoidal function of the duration allocated to the task.  Each task in the queue loses its value continuously. We elucidate on this trade-off and determine optimal policies for the human operator. We show the optimal policy requires the human to drop some tasks. We present a receding horizon optimization strategy, and compare it with the greedy policy. 
\end{abstract}


\section{Introduction}\label{sec:introduction}
The advent of cheap camera sensors, has prompted the extensive deployment of camera sensor networks for surveillance~\cite{WMB:09,CD:10}. Typically, the feeds of these cameras are send to a central location, where a human operator looks at these feeds and decides on the existence of certain features. The correctness of the decision on a particular feed depends on the time-duration the human operator allocates to it.  Thus, for a busy human operator, the optimal time-duration allocation is critically important. 

Recently, there has been significant interest in understanding the physics of  human decision making~\cite{RB-EB-etal:06}. Several mathematical models for human decision making have been proposed~\cite{RB-EB-etal:06, RWP:69}.  These models suggest that the correctness of the decision of a human operator in a binary decision making scenario evolves as a sigmoidal function of the time-duration allocated for decision. Thus, the probability of the correct decision by a human operator remains small till a critical time, and then jumps to a larger value.  When a human operator has to serve a queue of decision making tasks, then the tasks (e.g., feeds from camera) waiting in the queue lose value  continuously.  This trade-off between the correctness of the decision and the loss in the value of the pending tasks is of critical importance for the performance of the human operator.  In this paper, we address this trade-off, and determine optimal duration allocation policies for the human operator serving a decision making queue. Alternatively, we determine the task release rate that yields the desired accuracy for each task.

There has been a significant interest in the study of the performance of a human operator serving a queue.  Schmidt~\cite{DKS:78} model the human as a server and numerically study a queueing model to determine the performance of a human air traffic controller.  Maximally stabilizing task release policies for a human in the loop queue has been studied by Savla et al~\cite{KS-EF:10b, KS-EF:10, KS-EF:10a}.  Bertuccelli et al~\cite{LFB-NP-MLC:10} and Savla et al~\cite{KS-TT-EF:08} study the human supervisory control for the unmanned aerial vehicle operations.  
Donohue et al~\cite{MD-CL:09} study an optimal duration allocation  problem  is a time constrained static queue.

The optimal control of queueing systems~\cite{LIS:99} is a classical problem in queueing theory.  Stidham et al~\cite{SSJ-RRW:89} study the optimal service  policies  for a M/G/1 queue. They formulate a semi-Markov decision process, and describe the qualitative features of the solution. Certain assumptions in~\cite{SSJ-RRW:89} are relaxed by George et al~\cite{JMG-JMH:01}. Hern{\'a}dez-Lerma et al~\cite{OHL-SIM:83} study an optimal adaptive service rate control policy for a M/G/1 queue with an unknown arrival rate. 

We study the problem of optimal time-duration  allocation in a queue of binary decision making tasks with a human operator. We refer to such queues as {\it decision making queues}.  We consider three particular problems. First, we consider a time constrained static queue, where the human operator has to perform a given number of tasks with in a prescribed time. Second, we consider a static queue with latency penalty. Here, the human operator has to perform a given number of tasks. The operator incurs a penalty due to the delay in processing of each task. This penalty can be thought of as the loss in the value of the task over time.
Last, we consider a dynamic queue of the decision making tasks. The tasks arrive at a fixed rate and the operator incurs a penalty for the delay in processing of each task. In all the three problems, there is a trade-off between the reward obtained by processing a task, and the penalty incurred due to resulting delay in processing other tasks. We address this particular trade-off.  
Major contributions of this work are:

\begin{enumerate}
\item We determine a closed form optimal duration allocation policy for the time constrained static decision making queue and the static decision making queue with penalty. 
\item We provide a simple procedure to determine an optimal duration allocation policy for the dynamic decision making queue.
\item We rigorously establish that the optimal duration allocation policy may drop some tasks, i.e, not process the tasks at all. 
\end{enumerate}

The remainder of the paper is organized in the following way. We discuss some preliminary concepts  in Section~\ref{sec:preliminaries}. We present the problem setup  in Section~\ref{sec:problem-setup}.  We present the optimal allocation policy for  time constrained static queue  in Section~\ref{sec:static-constrained}. The static queue with latency penalty  is considered in Section~\ref{sec:static-latency}.  We present the optimal allocation policy  for the dynamic queue with latency penalty in Section~\ref{sec:dynamic-penalty}. 
We elucidate on these problems further through some examples in 
Section~\ref{sec:numerical-examples}. Our conclusions are presented in Section~\ref{sec:conclusion}.


\section{Preliminaries}\label{sec:preliminaries}
\subsection{Speed-accuracy trade-off in human decision making}

Consider the scenario where the human has to decide on one of the two alternatives $H_0$ and $H_1$, based on the collected evidence. 
 The evolution of the probability of correct decision has been studied in cognitive psychology literature~\cite{KLB-RRT:86, RWP:69, RB-EB-etal:06}. 
\begin{LaTeXdescription}
\item[ Pew's model:] The probability of deciding on hypothesis $H_1$, given that hypothesis $H_1$ is true, at a given time $t\in\real_{\ge0}$ is given by
\begin{equation*}
\prob(\text{say }H_1|H_1,t)= \frac{p_0}{1+e^{-(at-b)}},
\end{equation*}
where $p_0\in[0,1]$, $a,b\in\real$ are some parameters which depend on the human operator~\cite{KLB-RRT:86, RWP:69}.
\item[Drift diffusion model:] Conditioned on the hypothesis $H_1$, the evolution of the evidence for decision is modeled as a 
drift-diffusion process~\cite{RB-EB-etal:06}. Given a drift rate $\beta>0$, and diffusion rate $\sigma$, with a decision threshold $\eta$, the conditional probability of the correct decision is
\[
\prob(\text{say }H_1|H_1,t) =\frac{1}{\sqrt{2\pi\sigma^2 t}} \int_{\eta}^{\infty} e^{\frac{-(\Lambda -\beta t)^2}{2\sigma^2 t} }d\Lambda,
\]
where $\Lambda \equiv \mathcal{N}(\beta t, \sigma^2 t)$ is the evidence at time $t$.
\end{LaTeXdescription}

\begin{figure}[ht]
 \centering
     \setlength{\fboxrule }{0pt}
        \setlength{\fboxsep}{4pt}
\subfigure[Pew's model]{\label{fig:pew-model} 
        \fbox{ \includegraphics[width=0.21\textwidth]{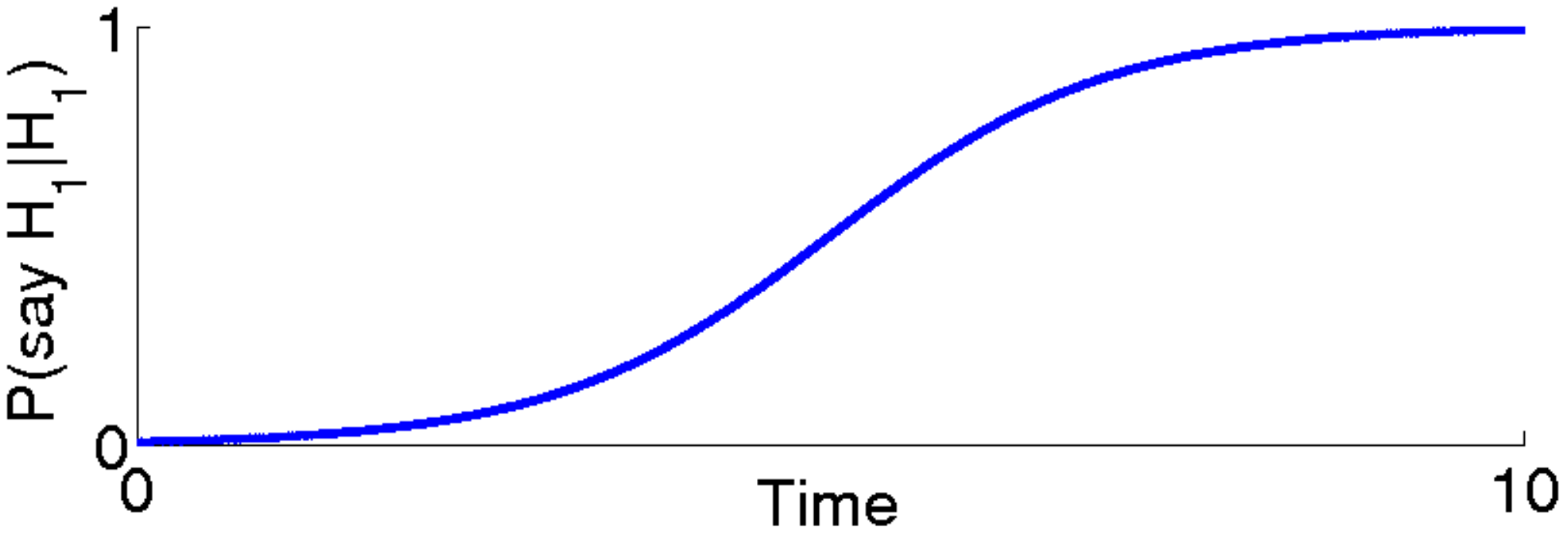}}}
\subfigure[Drift diffusion model]{\label{fig:ddm}\fbox{\includegraphics[width=0.21\textwidth]{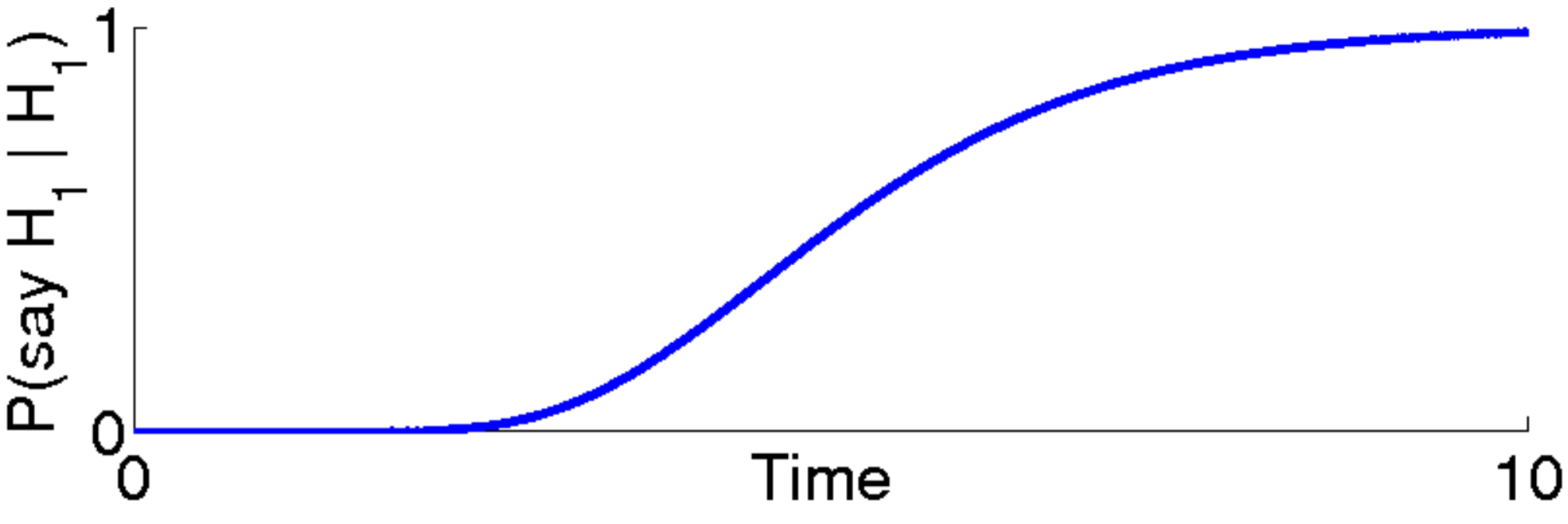}} }
\caption{The sigmoidal evolution of the probabilities of correct decision}
\end{figure} 

\subsection{Sigmoidal functions}
\begin{figure}[ht]
 \centering
     \setlength{\fboxrule }{0pt}
        \setlength{\fboxsep}{4pt}
\subfigure[]{\label{fig:sigmoidal-penalty}
        \fbox{ \includegraphics[width=0.21\textwidth]{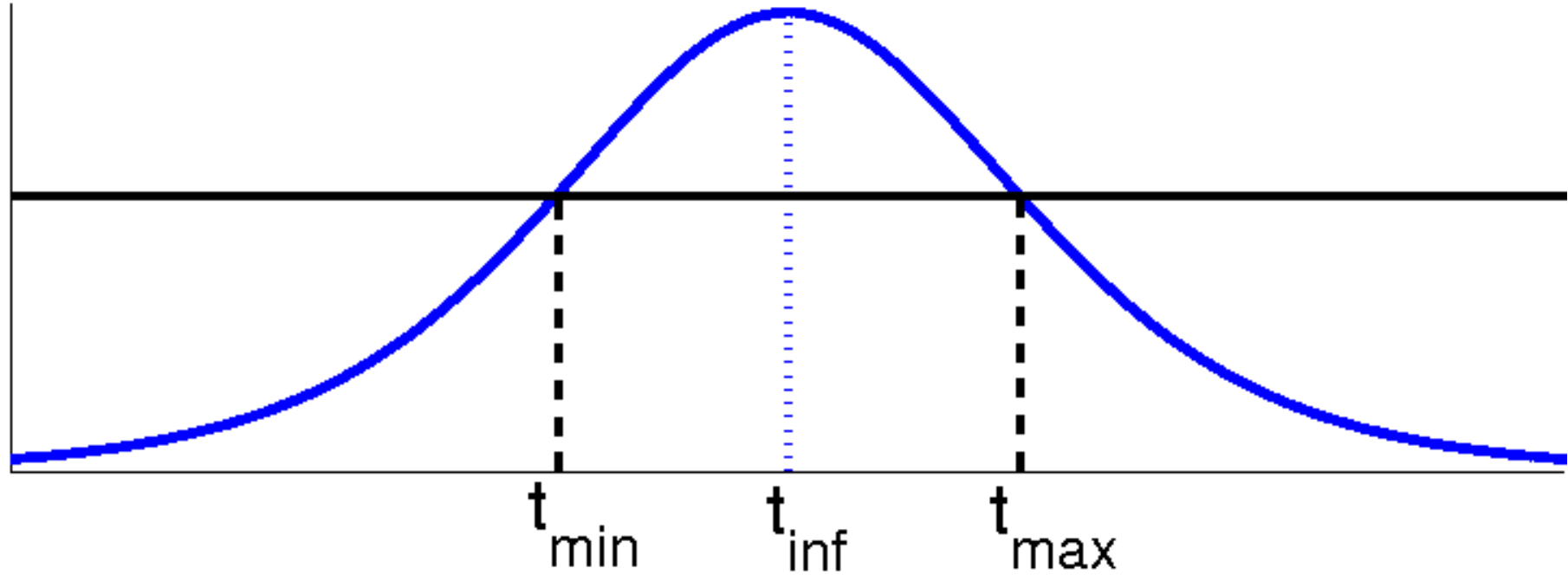}}}
\subfigure[]{\label{fig:sigmoidal_dder}\fbox{\includegraphics[width=0.21\textwidth]{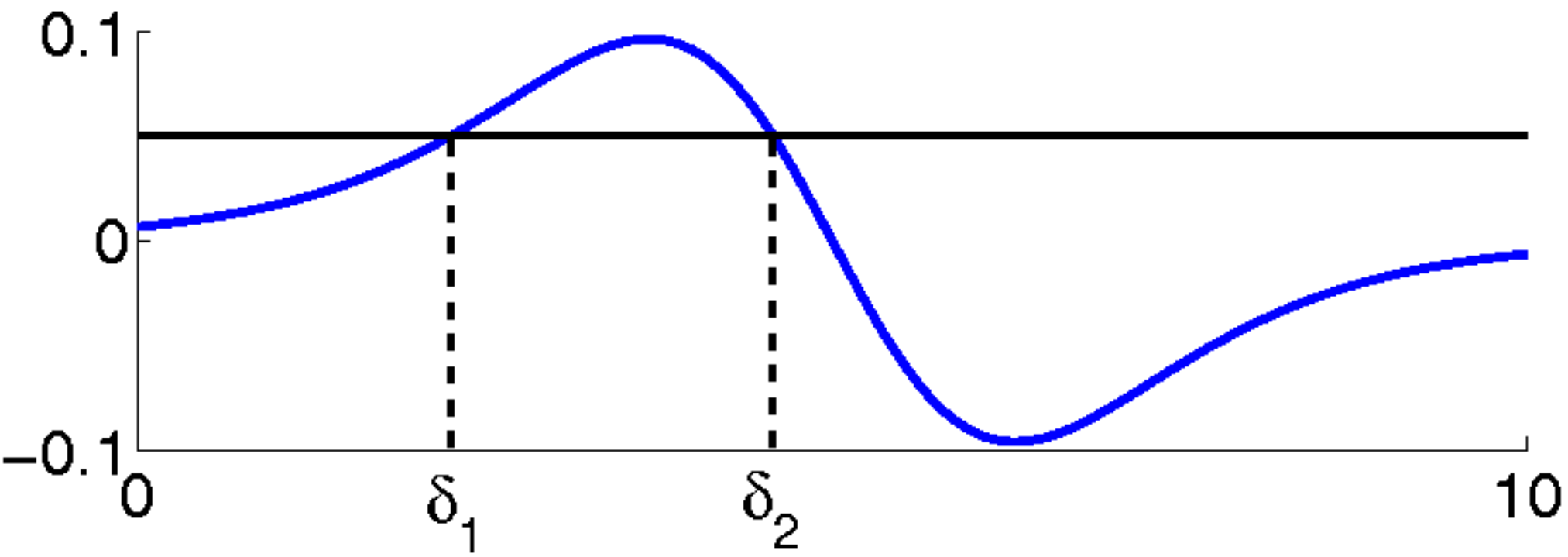}} }
\caption{(a) First derivative of the sigmoidal function and penalty rate. A particular value of the derivative may be attained at two different times. The total benefit decreases till $t_{\min}$, increases from $t_{\min}$ to $t_{\max}$, and then decreases again.  (b) Second derivative of the sigmoidal function. A particular positive value of the second derivative may be attained at two different times.}\label{fig:sigmoidal-derivatives}
\end{figure}

A smooth function $\map{f}{\real_{\geq0}}{\real_{\geq 0}}$ defined by
\[
f(t) = f_{\text{cvx}}(t) \mathcal{I}(t<t_{\text{inf}}) +f_{\text{cnv}}(t) \mathcal{I}(t\ge t_{\text{inf}}),
\]
where $ f_{\text{cvx}}$ and $f_{\text{cnv}}$ are  monotonically increasing  convex and concave functions, respectively, and $t_{\text{inf}}$ is the
inflection point. Derivative of sigmoidal function is unimodal with maximum at $t_{\text{inf}}$. Further, $f'(0)\ge 0$ and $\lim_{t\to \infty} f'(t)=0.$ Also, $\lim_{t\to \infty} f''(t)=0.$ A typical graph of the first and second derivative of a sigmoidal function is shown in Figure~\ref{fig:sigmoidal-derivatives}.
 Note that the evolution of the conditional 
probabilities of correct decision are sigmoidal functions in Pew's as well as drift-diffusion model.

\subsection{Receding horizon optimization}
Consider the following infinite horizon dynamic optimization problem.
\begin{equation}\label{eq:infinte-horizon-optimization}
\begin{split}
\maximize &\quad \sum_{\ell=1}^{\infty} \psi(x(\ell), u(\ell))\\
\subject &\quad x(\ell+1)= \phi(x(\ell), u(\ell)),
\end{split}
\end{equation}
where $x(\ell),u(\ell)\in\real$ are the state and control input at time $\ell\in\naturals$, respectively, $\map{\psi}{\real\times\real}{\real}$ is the stage cost, and 
$\map{\phi}{\real\times\real}{\real}$ defines the nonlinear evolution of the system.

The receding horizon optimization~\cite{EFC-CB:04} approximates the optimization problem~\eqref{eq:infinte-horizon-optimization} as the
 following finite horizon optimization problem at each 
stage $\theta\in\naturals$: 
\begin{equation}\label{eq:receding-horizon-optimization}
\begin{split}
\maximize &\quad \sum_{\ell=\theta}^{\theta+N} \psi(x(\ell), u(\ell))\\
\subject &\quad x(\ell+1)= \phi(x(\ell), u(\ell)),
\end{split}
\end{equation}
where $N\in\naturals$ is a finite horizon length. The receding horizon optimization is summarized as following
\begin{algorithm}

\caption{Receding horizon optimization}
\begin{algorithmic}[1]
  \STATE at time $\theta\in\natural$, observe state $x(\theta)$
  \STATE Solve optimal control problem~\eqref{eq:receding-horizon-optimization} and compute the optimal control inputs $u^*(\theta),\ldots,u^*(\theta+T)$\\
  \STATE Apply $u^*(\theta)$, and set $\theta=\theta+1$\\
  \STATE Go to step~\algostep{1}
\end{algorithmic}
\label{algo:receding-horizon}
\end{algorithm}

\section{Problem setup}\label{sec:problem-setup}
\begin{figure}[ht]
 \centering
\includegraphics[width=0.4\textwidth]{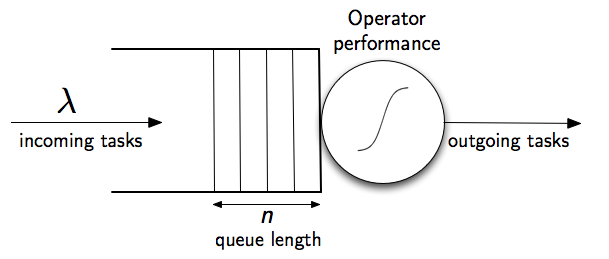}
\caption{Problem setup. The decision making tasks arrive at a rate $\lambda$. These tasks are served by a human operator with sigmoidal performance. Each task loses value while waiting in the queue. } \label{fig:problem-setup}
\end{figure}

We consider the problem of optimal time duration  allocation for a human operator. The decision making tasks arrive at a given rate and are processed by a human operator (see Figure~\ref{fig:problem-setup}.)  
The human receives a unit reward for the correct decision, while there is no penalty for a wrong decision.  For a decision made at time $t$, the expected reward is
\begin{equation}\label{eq:expected-reward}
\expt[\ones_{\text{say }H_1|H_1,t}] =\prob(\text{say }H_1|H_1,t)=f(t),
\end{equation}   
where $\map{f}{\real_{\ge0}}{[0,1]}$ is a sigmoidal function.

We consider three particular problems. First, we consider a time constrained static queue, i.e.,  
the human operator has to perform $N\in\naturals$ identical decision making 
tasks with in  time $T\in\real_{>0}$.  Second, we consider a static queue with penalty, i.e., 
 the human operator has to perform $N\in\naturals$
 identical decision making tasks, but each task  lose value at a constant rate per unit delay in its processing. Last,  
 we consider a dynamic queue of identical decision making tasks where each task lose value at a constant rate per unit delay in its processing.  For such a decision making queue, we are interested in the optimal time-duration allocation to each task. Alternatively, we are interested in the task release rate that will result in the desired accuracy for each task.  We intend to design a decision support system that tells the human operator the optimal time-duration allocation for each task.

\section{Time constrained static queue}\label{sec:static-constrained}

\subsection{Problem description}
Consider that the human operator has to perform $N\in\naturals$ identical tasks with in  time $T\in\real_{>0}$. Let the human operator allocates duration $t_{\ell}$ to the task $\ell\in\until{N}$. The following optimization problem encapsulates the objective of the human operator:
\begin{equation} \label{eq:maximize-no-arrival}
  \begin{split}
    \underset{\t}{\maximize} &\quad \sum_{\ell=1}^{N} f(t_{\ell})\\
    \subject  &\quad  \sum_{\ell=1}^{N} t_{\ell} \le T\\
              &\quad \t \succeq 0,
  \end{split}
\end{equation}
where $\t=\{t_1,\dots,t_N\}$ is  the
duration allocation vector.

\subsection{Optimal policy}
We define the Lagrangian $\map{L}{\real^N_{>0}\times\real_{\ge 0} \times\real^N_{\ge0} }{\real}$ for the constrained optimization problem~\eqref{eq:maximize-no-arrival} by
\begin{equation*}
L(\t,\alpha,\m) = \sum_{\ell=1}^{N} f(t_{\ell}) + \alpha (T- \sum_{\ell=1}^N t_{\ell}) + \m^{\trans} \t.
\end{equation*}
We also define the dual cost function $\map{h}{\real_{\ge 0}\times\real^N_{\ge0}}{\real}$ by 
\begin{equation}\label{eq:dual-cost}
h(\alpha,\m)= \underset{\t\in\real^N_{\ge 0}}\max L(\t,\alpha,\m).
\end{equation}
The dual problem to the primal optimization problem~\eqref{eq:maximize-no-arrival} is 
\begin{align}\label{eq:dual-no-arrival}
\begin{split}
\minimize &\quad h(\alpha,\m)\\ 
\subject &\quad \alpha \ge 0, \; \m\succeq 0.
\end{split}
\end{align}

\begin{lemma}[Constraint qualification]\label{lem:strong-duality}
For the primal problem~\eqref{eq:maximize-no-arrival} and the dual problem~\eqref{eq:dual-no-arrival}, the constraint qualification holds.
\end{lemma} \smallskip
\begin{proof}
Let $\t^*$ be an optimal solution to the primal problem~\eqref{eq:maximize-no-arrival}. 
 Without loss of generality assume that $t^*_1>0$, i.e., $t^*_1\ge 0$ is not an active constraint. 
 The gradients of the other constraints are $\e_j$, $j\in\{2,\ldots,N\}$ and $\ones_N$, where $\e_j$ is the ${j}^\text{th}$ 
 standard Cartesian unit vector, and $\ones_N$ is the $N$-tuple with all ones. It immediately follows that these gradients are linearly independent. 
 Thus, the constraint qualification~\cite{UD:08} holds.
\end{proof}

%

\begin{theorem}[Time constrained static queue]\label{thm:optimal-allocation-no-arrival}
For the optimization problem~\eqref{eq:maximize-no-arrival},
the optimal duration allocation vector $\t^*$ is an $N$-tuple with $m^*$ entries equal to $T/m^*$ and all other entries zero, where 
\begin{equation*}
m^*=\underset{m\in\until{N}}{\argmax} mf(T/m).
\end{equation*}
\end{theorem} \smallskip\smallskip
\begin{proof}
We apply the Karush-Kuhn-Tucker necessary conditions~\cite{SB-LV:04} for an optimal solution $\t^*,\alpha^*,\m^*$.\\
{\it Linear dependence of gradients}
\begin{multline}\label{eq:kkt-no-arrival-1}
\frac{\partial L}{\partial t^*_{\ell}}(\t^*,\alpha^*,\m^*) = f'(t^*_{\ell})-\alpha^*+\mu^*_{\ell} =0, \\
\forall \ell\in\until{N}.
\end{multline}
{\it Feasibility of the solution}
\begin{align}\label{eq:time-constraint}
T-\ones_N^{\trans}\t^* &\ge 0.\\
\t^*&\succeq 0.
\end{align}
{\it Complementarity conditions}
\begin{align}\label{eq:compl-no-arrival-1}
\alpha^*(T-\ones_N^{\trans}\t^*) &= 0.\\\label{eq:compl-no-arrival-2}
\mu_{\ell}^*t^*_{\ell}&= 0, \; \forall \ell\in\until{N}.
\end{align}
{\it Non-negativity of the multipliers}
\begin{equation}
\alpha^*\ge 0, \quad \m^*\succeq 0. 
\end{equation}
Since $f$ is a strictly increasing function, the constraint~\eqref{eq:time-constraint} should be active, and thus, from complementarity condition~\eqref{eq:compl-no-arrival-1} 
$\alpha^*>0$. Further, from equation~\eqref{eq:compl-no-arrival-2}, if $t^*_{\ell}\ne 0$, then $\mu_{\ell}^*=0$. 
Assume that $t^*_{\ell}\ne 0$, then, for each $\ell\in \until{N}$, the equation~\eqref{eq:kkt-no-arrival-1} yields
\begin{align}\label{eq:first-order-cond}
f'(t^*_{\ell})=\alpha^*.
\end{align}
If the equation~\eqref{eq:first-order-cond} has no solution, then the
Lagrangian is decreasing function of $t^*_{\ell}$, and the optimal duration
allocation is zero.  But $t^*_{\ell}\ne 0$, by assumption. Thus, there
exists a solution of equation~\eqref{eq:first-order-cond}.  Since, the
derivative of a sigmoidal function is uni-modal, there exists two solutions
of equation~\eqref{eq:first-order-cond}.  We refer to these values of
$t^*_{\ell}$ as $t^-$ and $t^+$, with the understanding that $t^+\ge
t^-$. From Figure~\ref{fig:sigmoidal-penalty}, it follows that a local
minima exists at $t^-$, while a local maxima exists at $t^+$.  Further note
that the Lagrangian is a decoupled function of $t_{\ell},\;
\ell\in\until{N}$. Thus, the optimal $\t^*$ will contain no $t^-$ entry.
Moreover, each entry will be either zero or $t^+$. Let $m^*\le N$ entries
of optimal $\t^*$ be non-zero.  Any allocation with $m$ identical non-zero
entries, and remaining zero entries yields a reward $m f(T/m)$. Therefore,
the optimal non-zero entries are as stated in the theorem.
\end{proof}

\begin{remark}[Notes on concavity I]
If the performance function  $f$ is a concave function, the optimal policy for the time constrained queue is to allocate equal time to each task.\oprocend
\end{remark}

\section{Static queue with latency penalty}\label{sec:static-latency}

\subsection{Problem description}
Consider that the human operator has to perform $N\in\naturals$ identical tasks. Let the human operator assigns time $t_{\ell}$ to the task $\ell\in\until{N}$. 
 The operator receives an expected reward $f(t_{\ell})$ for an assignment $t_{\ell}$ to the task $\ell$, while she incurs a latency cost $c$ 
 per unit time for the delay in processing of each task.  The following optimization problem encapsulates the objective of the human operator:
\begin{equation} \label{eq:maximize-no-arrival-latency}
   \underset{\t\succeq0}{\maximize} \quad \frac{1}{N}\sum_{\ell=1}^{N} \big( f(t_{\ell})-ct_{\ell}(N-\ell +1)\big),\\
\end{equation}
where $\t=\{t_1,\dots,t_N\}$ is the duration allocation vector.

\subsection{Optimal policy}

\begin{theorem}[Static queue with latency penalty]\label{thm:optimal-allocation-no-arrival-penalty}
For the optimization problem~\eqref{eq:maximize-no-arrival-latency}, the optimal allocation at stage $\ell\in\until{N}$ is 
\begin{align*}
  t^*_{\ell} &:= \underset{\beta\in\{0,t^{\dag}_{\ell}\}}{\argmax} \big( f(\beta) -c (N\!-\ell\! +1)\beta\big), 
\end{align*}
where
\begin{align*}
     t^{\dag}_{\ell} := 
   \begin{cases} 
   0, & \quad\text{if }  f'(t_{\text{inf}}) < c(N\!-\!\ell\!+\!1),\\
   \max\{ t\in\real_{>0} \;|\;f'(t) = \!\!&\!\! c(N\!-\ell\! +1) \},\quad\text{otherwise.}
      \end{cases}
  \end{align*}
\end{theorem} \smallskip \smallskip
\begin{proof}
  The proof is similar to the proof of
  Theorem~\ref{thm:optimal-allocation-no-arrival}. If $f'(t_{\text{inf}}) <
  c(N\!-\!\ell\!+\!1)$, then the objective is a decreasing function of
  $t_{\ell}$, and the optimal allocation is zero. Otherwise, the
  local maximum  lies at the  intersection of  penalty rate $c(N- \ell+1)$ with the decreasing portion of the $f'$ (see Figure~\ref{fig:sigmoidal-penalty}).  The optimal allocation
  $t^*_{\ell}$ is determined by comparing value of the objective function
  at the local maxima and the boundary.
\end{proof}

\begin{remark}[Notes on concavity II]
The optimal duration allocation for the static queue with latency penalty decreases to a critical value $t_{\text{crit}}>\tinf$ with increasing penalty rate, then jumps down to zero. Instead, if  the performance function $f$ is concave, then the optimal duration allocation decreases continuously to zero with increasing penalty rate.  \oprocend
\end{remark}

\section{Dynamic queue with latency penalty }\label{sec:dynamic-penalty}

\subsection{Problem description}
Consider that the human operator has to serve  a queue of identical decision making tasks. Assume that the tasks arrive as a 
Poisson's process with rate $\lambda>0$. We define the processing of job $\ell\in\naturals$ as the stage $\ell$. 
Let the operator assign deterministic time $t_{\ell}$ at stage $\ell$.  The operator receives an expected reward $f(t_{\ell})$ 
for a duration allocation $t_{\ell}$ at stage $\ell$, while she incurs a latency cost $c$ per unit time for the delay in processing each task.
 We denote the queue length at stage $\ell$ by $n_{\ell}$. The objective of the operator is to maximize her  infinite horizon expected reward. 
  The following optimization problem encapsulates the objective of the human operator:
\begin{equation}\label{eq:maximize-with-arrivals}
\begin{split}
{\maximize}& \quad \lim_{L\to \infty}\frac{1}{L}\sum_{\ell=1}^{L}
\Big(f(t_{\ell}) - c\expt[n_{\ell}]t_{\ell}-\frac{c\lambda
  t_{\ell}^2}{2}\Big)\\ 
\subject&\quad \expt[n_{\ell+1}]= \max\{0, \expt[n_{\ell}] -1 +\lambda t_{\ell} \}\\
&\quad t_{\ell} \ge 0, \forall \ell\in \naturals,  
\end{split}
\end{equation}
 where the term $c\lambda t_{\ell}^2/2$ is the expected penalty due to the tasks arriving during stage $\ell$.


\subsection{Optimal policy}
\smallskip\noindent
\textit{Finite horizon optimization}
\smallskip

We wish to approximately solve the infinite horizon optimization
problem~\eqref{eq:maximize-with-arrivals} using receding horizon optimization as
in Algorithm~\ref{algo:receding-horizon}.  To do so, we pose the following
finite horizon optimization problem:
\begin{equation}\label{eq:maximize-receding-horizon}
\begin{split}
\underset{\t\succeq0}{\maximize}& \quad \frac{1}{N}\sum_{\ell=1}^{N} \Big(f(t_{\ell}) - c\expt[n_{\ell}]t_{\ell}-\frac{c\lambda t_{\ell}^2}{2}\Big)\\
\subject&\quad \expt[n_{\ell+1}]= \max\{0, \expt[n_{\ell}] -1 +\lambda t_{\ell} \},
\end{split}
\end{equation}
where $\t=\{t_1,\dots,t_N\}$ is the duration allocation vector.

Without loss of generality, we assume that the queue length is always
non-zero. If the queue is empty at some time, then the operator will wait till the arrival of new task. There is no explicit penalty for the operator to be idle.
Under this assumption we have
\[
\expt[n_{\ell}] =n_1-\ell +1 +\lambda \sum_{j=1}^{\ell-1} t_j.
\]
Some algebraic manipulations show that the objective function of the
optimization problem~\eqref{eq:maximize-receding-horizon} is equivalent to
the function $\map{J}{\real_{\ge 0}^N}{\real}$ defined by
\[
J(\t):=\frac{1}{N}\sum_{\ell=1}^{N} \Big( f(t_{\ell})
-c(n_1-\ell+1)t_{\ell}-c\lambda t_{\ell} \sum_{j=1}^N t_j + \frac{c \lambda
  t_{\ell}^2}{2} \Big),
\]
where  $c$ is the
penalty rate, $\lambda$ is the arrival rate, and $n_1$ is the initial queue
length.  Thus, the optimization
problem~\eqref{eq:maximize-receding-horizon} is equivalent to
\begin{equation}\label{eq:optimize-equivalent}
\underset{\t\succeq0}{\maximize}\quad  J(\t).
\end{equation}



We define $\map{\fdi}{[0, f'(\tinf)]}{\real_{\ge 0}}$ by
\[
\fdi(y) = \max\setdef{t\in\real_{\geq0}}{f(t)=y}. 
\]
Note that the definition of $\fdi$ is consistent with Figure~\ref{fig:sigmoidal-penalty}.

For the optimization problem~\eqref{eq:maximize-receding-horizon}, assume
that the optimal policy allocates a strictly positive time only to the
tasks in the set $\Tproc\subseteq \until{N}$, which we call the \emph{set
  of processed tasks}.  (Accordingly, the policy allocates zero time to the
tasks in $\until{N}\setminus\Tproc$). Without loss of generality, assume 
\[
\Tproc:=\{\eta_1,\ldots,\eta_m\},
\] 
where $\eta_1<\cdots<\eta_m$ and $m\leq N$.  A duration allocation vector $\t$ is said to be consistent with  $\Tproc$ if  only the tasks in  $\Tproc$ are allocated non-zero duration.

\begin{figure}[ht]
 \centering
     \setlength{\fboxrule }{0pt}
        \setlength{\fboxsep}{4pt}
\subfigure[]{\label{fig:inverse-sigmoidal} 
        \fbox{\includegraphics[width=0.21\textwidth]{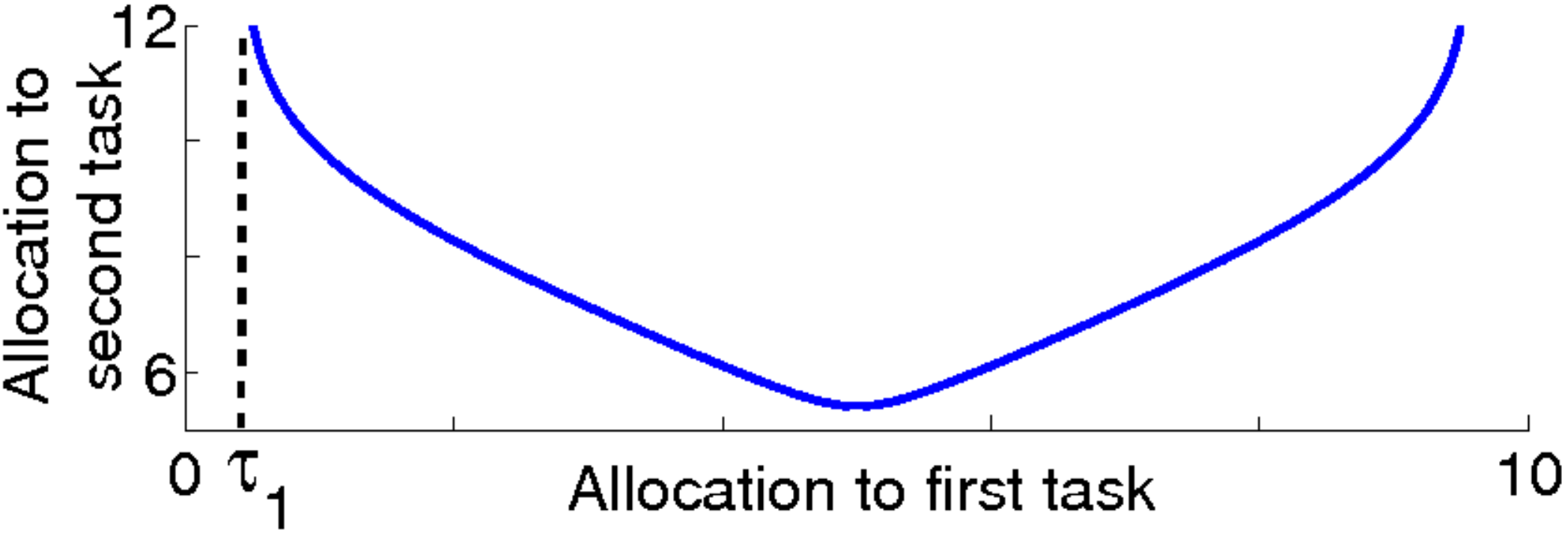}}}
\subfigure[]{\label{fig:reward-penalty}\fbox{\includegraphics[width=0.21\textwidth]{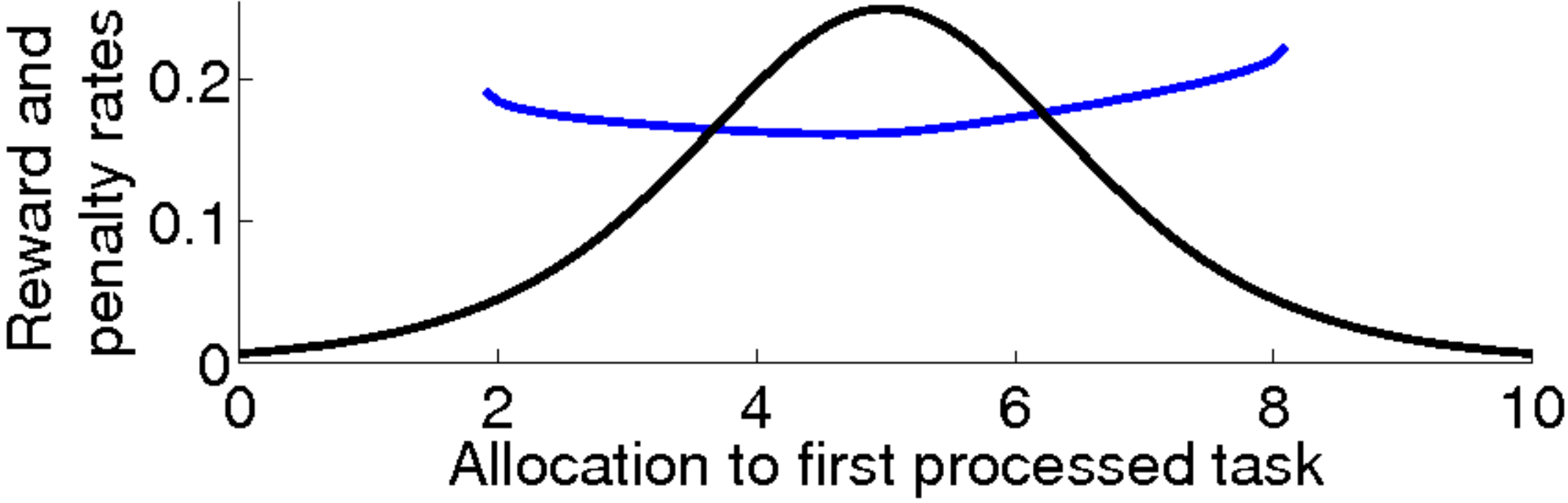}} }
\caption{(a) Feasible allocations to second task parametrized by allocation to first task. (b) The penalty rate and the reward  rate due to the allocation to the first task.}
\end{figure}

\begin{lemma}[Properties of  maximum point]\label{lem:penalty} 

 For the optimization problem~\eqref{eq:optimize-equivalent}, and  a set of processed tasks $\Tproc$, the following statements hold:
 \begin{enumerate}
 \item  A global maximum point $\t^*$ satisfy $t^*_{\eta_j} \ge t^*_{\eta_k}$, for $j\ge k$, $j,k\in\until{m}$. 
 \item   A  local maximum point $\t^{\dag}$ consistent with $\Tproc$ satisfies
 \begin{equation}\label{eq:objective-derivative}
 f'(t_{\eta_k}^{\dag})= c(n_1-\eta_k+1) + c\lambda \sum_{i=1}^{m} t^{\dag}_{\eta_i}, \forall k\in\until{m}.
 \end{equation}
 \item The system on equations~\eqref{eq:objective-derivative}, can be reduced to 
 \[
 f'(t_{\eta_1}^\dag) = \mathcal{P}(t_{\eta_1}^\dag),\text{ and } t^\dag_{\eta_k} =\fdi(f'(t_{\eta_1}^\dag)-c(\eta_k-\eta_1)),
 \]
 for each $k\in\{2,\dots,m\}$, where
 $\map{\mathcal{P}}{\real_{> 0}}{\real\union\{+\infty\}}$ is defined by
\[
\mathcal{P}(t)=\begin{cases}
p(t), \quad & \text{if }  f'(t)\ge c(\eta_m-\eta_1),\\
+\infty,& \text{otherwise,}
\end{cases}
\]
where \\
$\displaystyle p(t)\!=\!c\Big(\!n_1-\eta_1+1+  \lambda t+\lambda\!\sum_{k=2}^m\!
\fdi\!\big(f'(t) -c(\eta_k-\eta_1)\!\big)\! \Big)$.
\item  A local maximum point $\t^{\dag}$ consistent with $\Tproc$ satisfies
\[  
  f''(t_{\eta_k})\le c\lambda, \text{ for each } k\in\until{m}.
 \] 
 \end{enumerate}
 \end{lemma} \smallskip
 \begin{proof}
 We start by proving the first statement.   Assume $t_{\eta_j}^*<t_{\eta_k}^*$ and define the allocation vector $\bar{\t}$ consistent with $\Tproc$ by
  \[
  \bar{t}_{\eta_i}=\begin{cases}
    t^*_{\eta_i}, & \text{if } i\in\until{m}\setminus\{j,k\},\\
    t^*_{\eta_j},& \text{if } i=k,\\
    t^*_{\eta_k},&\text{if } i=j.
  \end{cases}
  \]
  It is easy to see that
  \[
  J(\t^*)-J(\bar{\t})  = (\eta_j-\eta_k)(t_{\eta_j}^*-t_{\eta_k}^*)<0.
  \]
  This inequality contradicts the assumption that $\t^*$ is a solution of
  the optimization problem~\eqref{eq:optimize-equivalent}.
  
   To prove the second statement, note that a local maximum is achieved at the boundary of the feasible region or at
  the set where the Jacobian of $J$ is zero. At the boundary of the
  feasible region $\real^N_{\ge 0}$, some of the allocations are zero.   Given the $m$ non-zero allocations, the Jacobian of the function
  $J$ projected on the space spanned by the non-zero allocations must be zero.  The expressions in the
  theorem are obtained by setting the Jacobian to zero.
  
  To prove the third statement, we
  subtract, the expression in equation~\eqref{eq:objective-derivative} for $k=j$ from the expression for $k=1$ to get
\begin{align}\label{eq:allocation-contour}
 f'(t_{\eta_j}) = f'(t_{\eta_1}) - c(\eta_j-\eta_1).
\end{align}
There exists a solution of  equation~\eqref{eq:allocation-contour} if and only if $f'(t_{\eta_1})\ge c(\eta_j-\eta_1)$. If $f'(t_{\eta_1})< c(\eta_j-\eta_1)+f'(0)$, then there exists only one solution. Otherwise, there exist two solutions. It can be seen that if there exist two solutions $t^{\pm}_j$, $t_j^-<t_j^+$, then $t_j^-<t_{\eta_1}<t_j^+$.  From i), only possible allocation is $t_j^+$. Notice that $t_j^+=\fdi(f'(t_{\eta_1}) - c(\eta_j-\eta_1))$.  
This yields feasible time allocation to each task $\eta_j, j\in\{2,\ldots,m\}$ parametrized  by the time allocation to the task $\eta_1$. A typical allocation is shown in Figure~\ref{fig:inverse-sigmoidal}.
We further note  that the effective penalty rate for the task $\eta_1$ is $c(n_1-\eta_1+1)+c\lambda\sum_{j=1}^m t_{\eta_j}$.
Using the expression of $t_{\eta_j},j\in\{2,\ldots,m\}$, parametrized by $t_{\eta_1}$, we obtain the expression for  $\mathcal{P}$. 

To prove the last statement, we observe that
the Hessian  of the function $J$ is 
\[
\frac{\partial^2 J}{\partial \t^2} = \text{diag}(f''(t_{\eta_1}),\ldots, f''(t_{\eta_m})) - c \lambda \ones_m \ones_m^T,
\]
where $\text{diag}(\cdot)$ represents a diagonal matrix with the argument as diagonal entries. 
For a local maximum to exist at non-zero duration allocations $\{t_{\eta_1},\dots,t_{\eta_m}\}$, the Hessian must be negative semidefinite.
Thus, a necessary condition for Hessian to be negative semidefinite is that diagonal entries are non-positive.
 \end{proof}

We refer to the function $\mathcal{P}$ by the \emph{effective penalty rate} for the first processed task. A typical graph of  $\mathcal{P}$ is shown in Figure~\ref{fig:reward-penalty}.
Given $\Tproc$, a feasible allocation to the task $\eta_1$ is such that $f(t_{\eta_1})-c(\eta_j-\eta_1)>0$, for each $j\in\{2,\dots,m\}$.
For a given $\Tproc$, we define the minimum feasible duration allocated to task $\eta_1$ (see Figure~\ref{fig:inverse-sigmoidal}) by
\begin{align*}
  \tau_1\!:=\!\begin{cases}
    \!\min\setdef{ t\!\in\!\real_{\!\ge\! 0\!}\!}{ \! f'\!(t)\!=\!c(\!\eta_m\!-\!\eta_1\!)\!}, & \!\text{if }\! f'\!(\!\tinf\!)\!\ge\! c(\!\eta_m\!-\!\eta_1\!),\\
    0, & \!\text{otherwise.}
  \end{cases}
  \end{align*}


Let $f''_{\text{max}}$ be the maximum value of $f''$. We now define the points at which the function $f''-c\lambda$ changes its sign (see Figure~\ref{fig:sigmoidal_dder}): 
 \begin{align*} 
\delta_1&:=\begin{cases}
 \min\setdef{t\in\real_{\ge 0}}{f''(t)=c\lambda}, & \text{if }  c\lambda\in[f''(0),f''_{\text{max}}],\\
 0,&\text{otherwise,}
 \end{cases}
 \\
 \delta_2 &:=\begin{cases}
 \max\setdef{t\in\real_{\ge 0}}{f''(t)=c\lambda}, & \text{if }  c\lambda \le f''_{\text{max}},\\
 0,&\text{otherwise.}
 \end{cases}
\end{align*} 


\begin{theorem}[Dynamic queue with latency penalty]
  \label{thm:allocation-maximum}
  For the optimization problem~\eqref{eq:optimize-equivalent},
  consider a set of processed tasks $\Tproc$. The following statements
  hold:
  \begin{enumerate}
  \item there exists a local maximum point consistent with $\Tproc$, if 
    \begin{align} 
  \label{cond:large-large}
  f'(\delta_2) &\ge  \mathcal{P}(\delta_2);
  \end{align}
\item  there exists a local maximum point consistent with $\Tproc$, if 
    \begin{align} 
\label{cond:small-large}
  f'(\tau_1) &\le \mathcal{P}(\tau_1), \; f'(\delta_1) \ge \mathcal{P}(\delta_1), \;\text{and } \delta_1\ge \tau_1;
  \end{align}
  \item if both
    conditions~\eqref{cond:large-large} and~\eqref{cond:small-large} are
    false, then
  there exist no local maximum point consistent with $\Tproc$.
\end{enumerate}
\end{theorem} \smallskip
\begin{proof}
A critical allocation to task $\eta_1$ is located at the intersection of the graph of the reward rate function $f'(t_{\eta_1})$ to the penalty functions $\mathcal{P}(t_{\eta_1})$. 
From Lemma~\ref{lem:penalty}, a necessary condition for existence of local maximum at a critical point is $f''(t_{\eta_1}) \le c\lambda$ . 
Further, for $t_{\eta_1}\in{]0,\delta_1]}\union {[\delta_2,\infty[}$, $f''(t_{\eta_1}) \le c\lambda$. It can be seen that if condition~\eqref{cond:large-large} holds, then the reward function $f'(t_{\eta_1})$  and the effective penalty function $\mathcal{P}(t_{\eta_1})$ intersect in the region ${[\delta_2,\infty[}$. Similarly,  condition~\eqref{cond:small-large} ensure the intersection of the graph of the reward function $f'(t_{\eta_1})$ with the effective penalty function $\mathcal{P}(t_{\eta_1})$ in the region ${]0,\delta_1]}$. 
\end{proof}


We now provide a procedure to provide solution to the optimization problem~\eqref{eq:maximize-receding-horizon}. Given a sequence of zero and 
non-zero allocations $\xi\in\{0,+\}^N$, we denote 
the corresponding critical allocation for maximum by $\t(\xi)$. The details of the procedure are shown in Algorithm~\ref{algo:decision-making-queue}.
\begin{algorithm}[h]
\caption{Optimal allocation for decision making queue}
\begin{algorithmic}[1]
  \STATE  given $n_1$, $N$, $c$, $\lambda$
  \STATE  $k:=0$; $\mathcal{A}:=\phi$;
  \STATE  {\bf for each} string $\xi\in\{0,+\}^N$
  \STATE  \quad set $\Tproc:=
  \setdef{i\in\until{N}}{\xi_i= +}$ 
  \STATE  \quad {\bf if} condition~\eqref{cond:large-large}~or~\eqref{cond:small-large} 
  \STATE \quad {\bf then} determine critical allocations \\ 
\hfill   for maximum  $t_{\eta_1}^{\dag}$ via bisection  algorithm
\STATE \quad   \qquad \;determine  allocations  $t_{\eta_j}^{\dag}, j\in\{2,\ldots,m\}$
\STATE \quad \qquad \;determine expected queue lengths\\
\hfill  $\expt[n_{\ell}], \ell\in\until{N}$
\STATE  \quad   \qquad \; {\bf if} $\expt[n_{\ell}]>0, \forall \ell\in\until{N}$
\STATE  \quad   \qquad \; {\bf then} $\mathcal{A}=\mathcal{A}\union \{\t^{\dag}(\xi)\}$
  \STATE optimal allocation $\t^*=\text{argmax}_{\t\in\mathcal{A}} J(\t)$
\end{algorithmic}
\label{algo:decision-making-queue}
\end{algorithm}
\begin{remark}[Notes on concavity III]
With the increasing penalty rate as well as the increasing arrival rate, the time duration allocation decreases to a critical value $t_{\text{crit}} > \tinf$  and then jumps down to zero, 
for the dynamic queue with latency penalty. Instead, if the performance function $f$ is concave, then the duration allocation decreases continuously to zero with increasing penalty rate as well as increasing arrival rate. \oprocend
\end{remark}



\section{Numerical examples}\label{sec:numerical-examples}

\begin{figure}[ht]
 \centering
     \setlength{\fboxrule }{0pt}
        \setlength{\fboxsep}{4pt}
\subfigure[Time constrained static queue]{\label{fig:time-constrained-allocation} 
        \fbox{\includegraphics[width=0.21\textwidth]{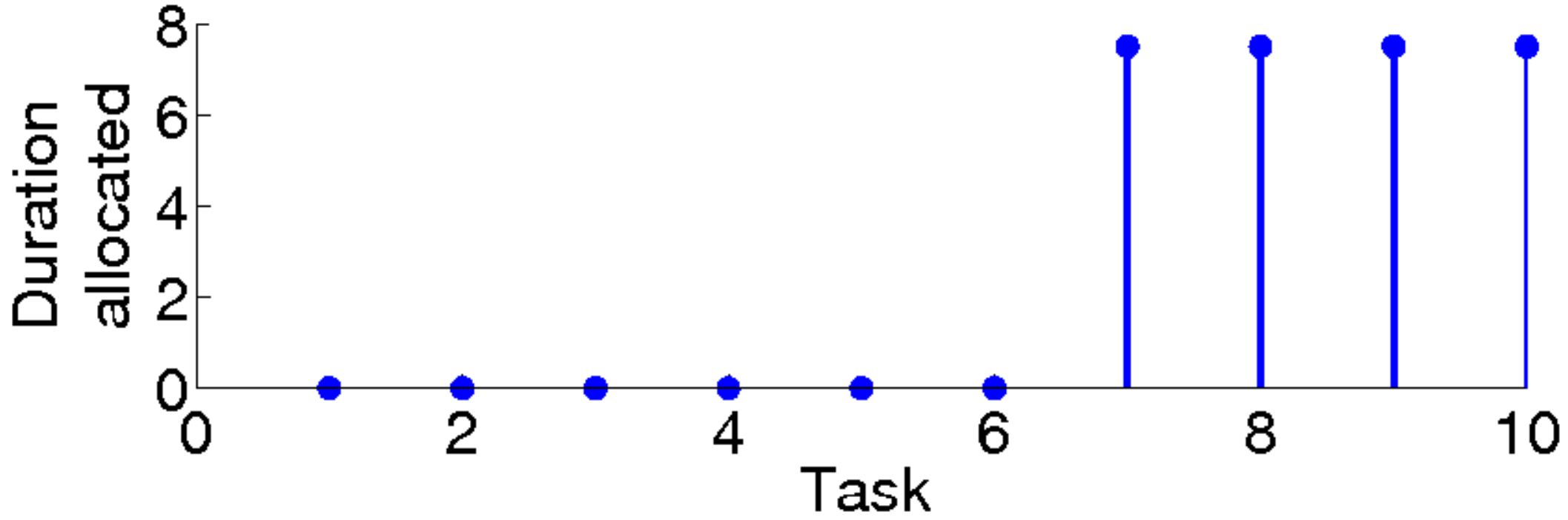}}}
\subfigure[Static queue wt. latency penalty]{\label{fig:static-allocation-penalty}\fbox{\includegraphics[width=0.21\textwidth]{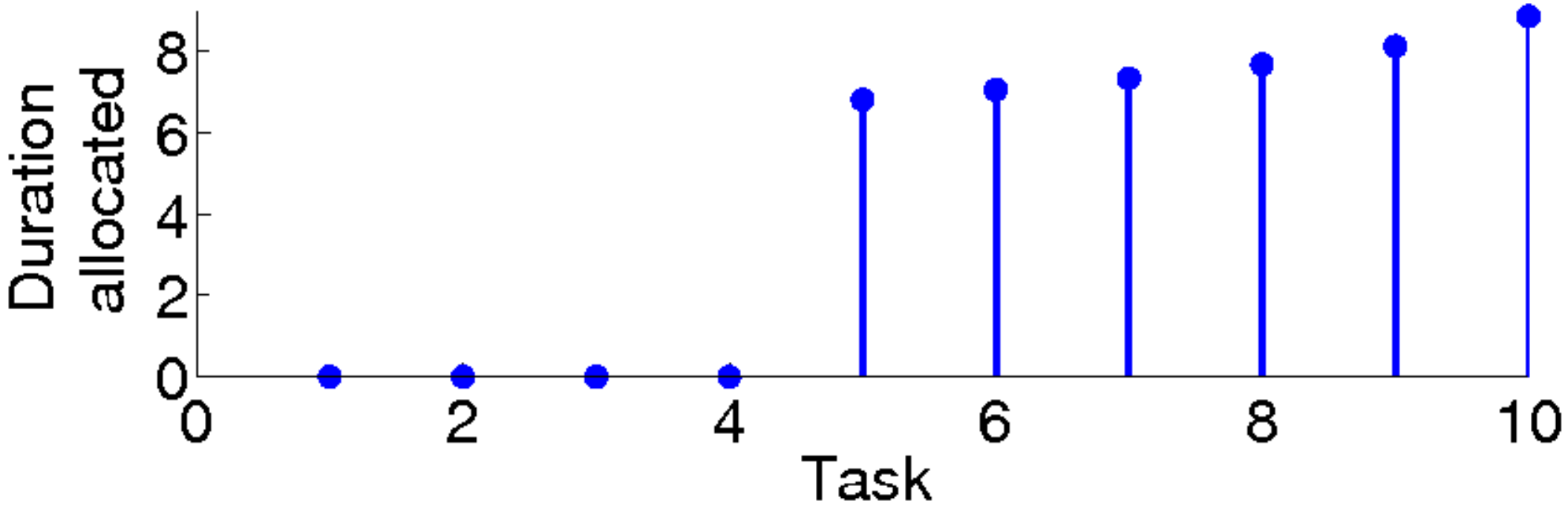}} }
\caption{Optimal allocations}
\end{figure}

We now elucidate on the optimal policies for the three problems through some numerical examples. We consider three examples.  In the first and second example, we 
demonstrate application of Theorem~\ref{thm:optimal-allocation-no-arrival}~and~\ref{thm:optimal-allocation-no-arrival-penalty}, respectively.  In the third example, the Algorithm~\ref{algo:decision-making-queue} is utilized in a receding horizon fashion. 
\begin{example}
If the human operator has to serve $N=10$ tasks in time $T=30$ secs, and the human receives an expected reward 
$f(t)=1/(1+\exp(5-t))$ for an allocation of duration $t$ secs to a task, then the 
optimal policy for the human is to drop six tasks, and allocate $7.5$ secs to any four tasks. 
An optimal allocation is shown in Figure~\ref{fig:time-constrained-allocation}.\oprocend
\end{example}

\begin{example}
If the human operator has to serve $N=10$ tasks and the human receives an expected reward $f(t)=1/(1+\exp(5-t))$
 for an allocation of duration $t$ secs to a task, while she incurs a penalty 
 $c=0.02$ per sec for each pending task,  then the optimal policy for the human is shown in 
 Figure~\ref{fig:static-allocation-penalty} .  Note that the optimal duration allocation increases with decreasing number of pending tasks.\oprocend
\end{example}

\begin{figure}[ht]
 \centering
     \setlength{\fboxrule }{0pt}
        \setlength{\fboxsep}{4pt}
\subfigure[Low arrival rate]{
        \fbox{\includegraphics[width=0.23\textwidth]{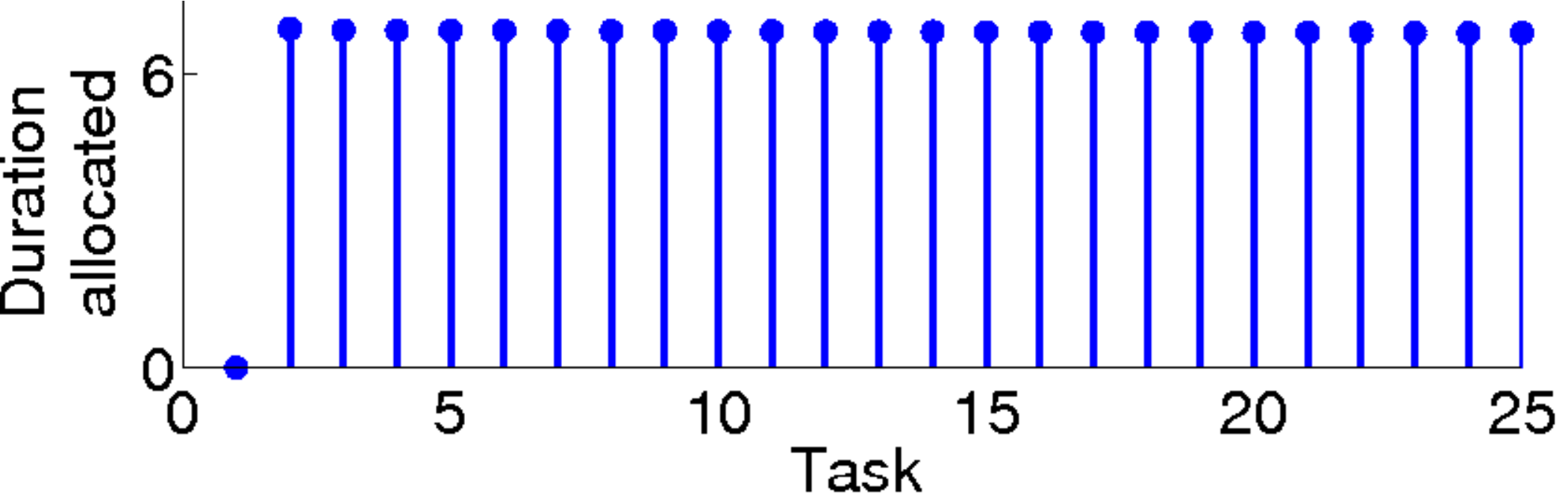}}
\fbox{\includegraphics[width=0.23\textwidth]{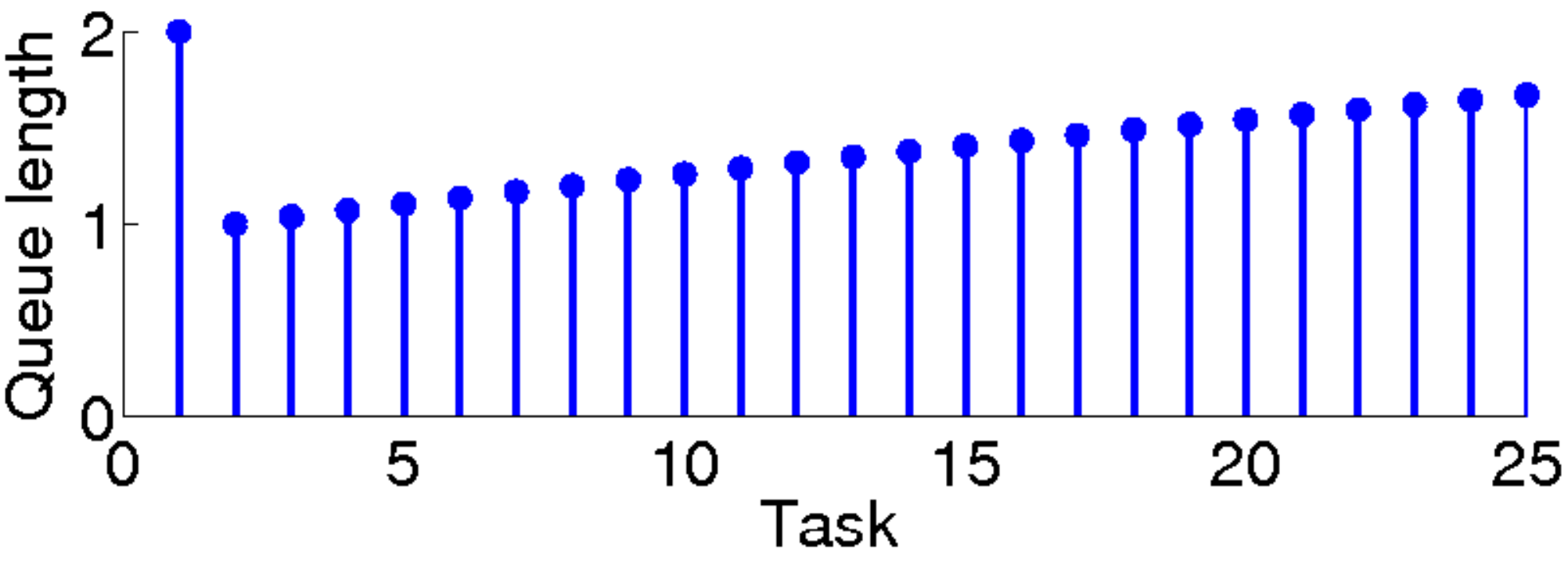}} }
\subfigure[Moderate arrival rate]{
        \fbox{\includegraphics[width=0.23\textwidth]{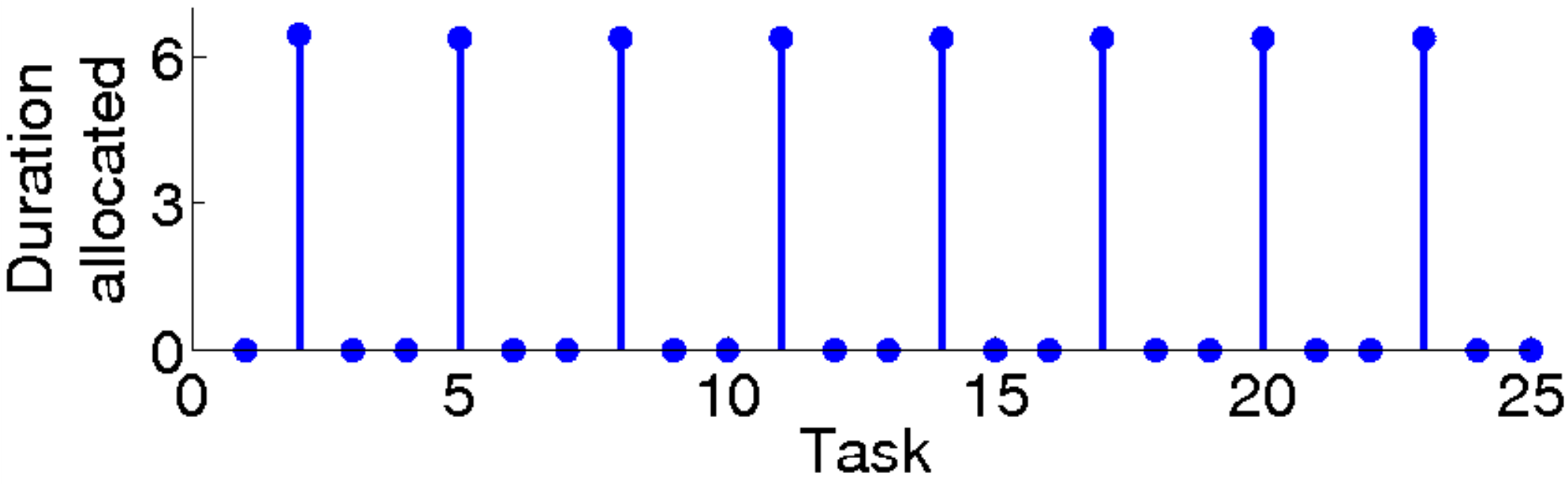}}
\fbox{\includegraphics[width=0.23\textwidth]{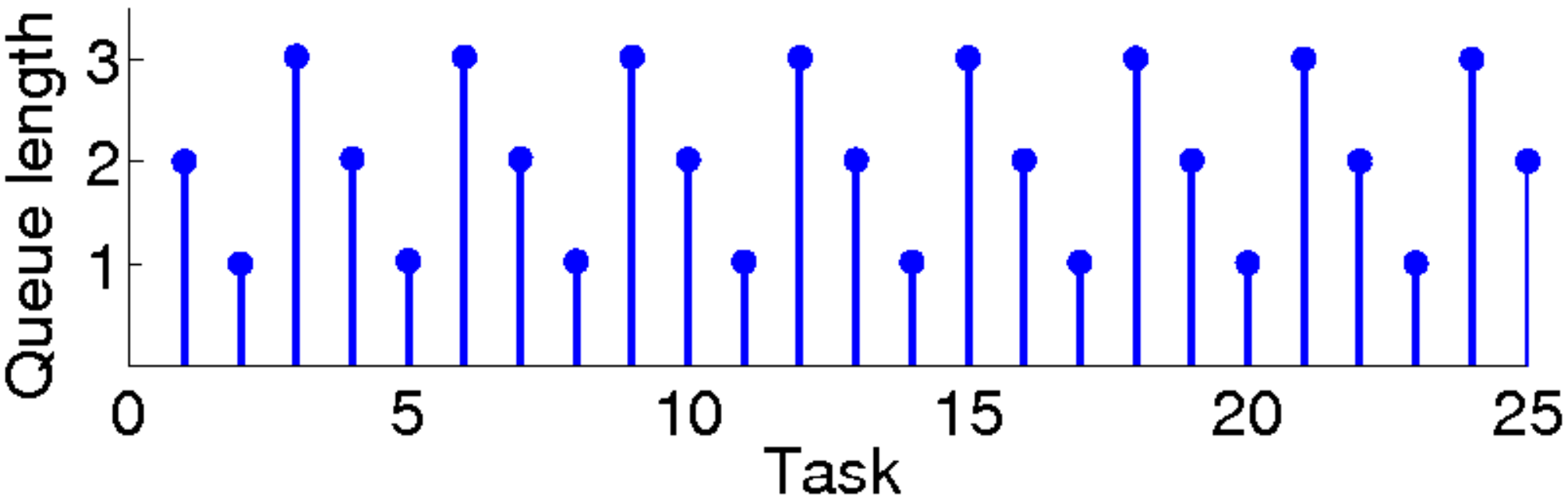}} }
\subfigure[High arrival rate]{
        \fbox{\includegraphics[width=0.23\textwidth]{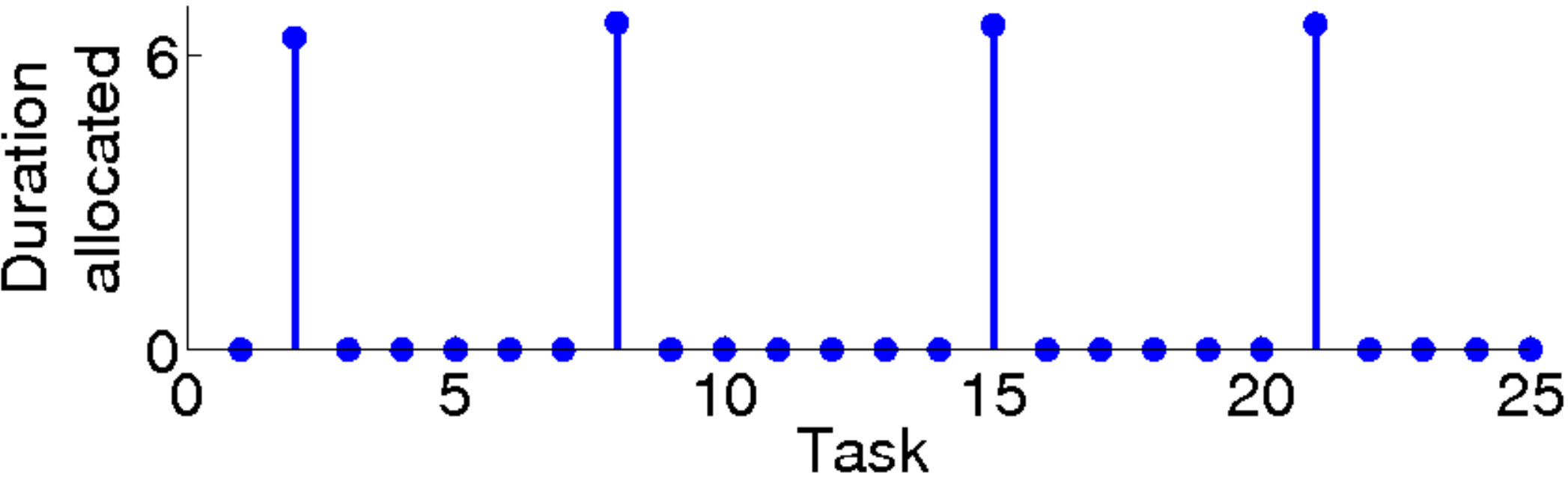}}
\fbox{\includegraphics[width=0.23\textwidth]{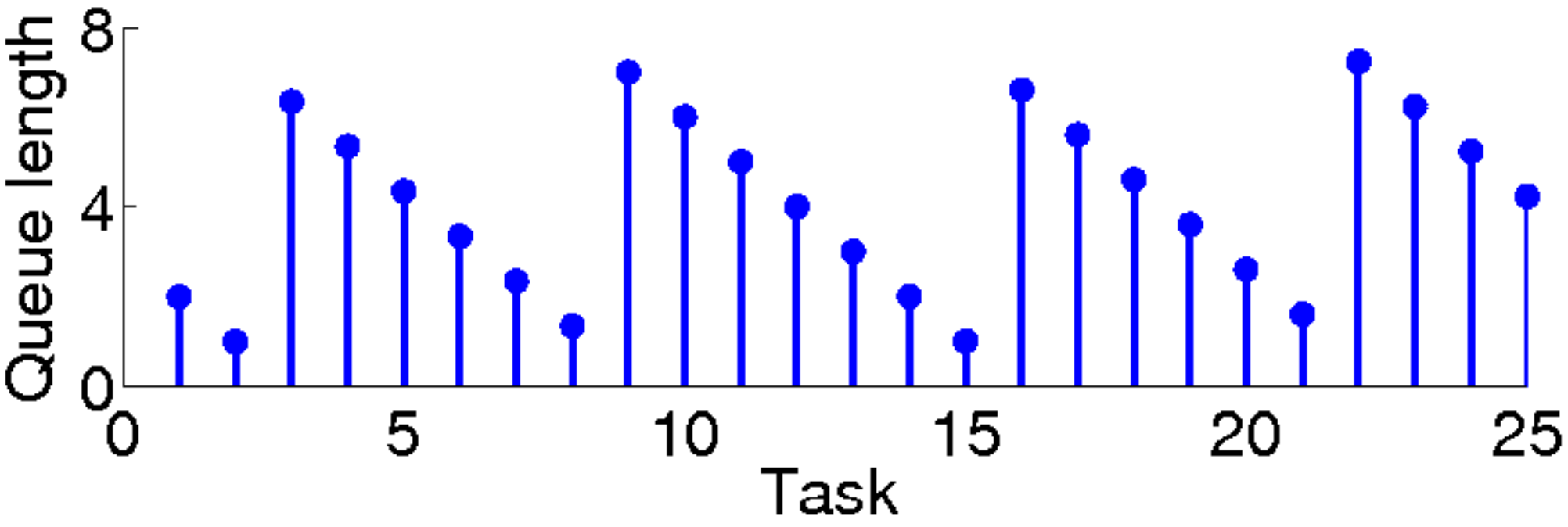}} }
\caption{Receding horizon policy for {\it expected evolution} of the dynamic queue with latency penalty. An optimization problem with horizon length $N=10$ is solved at each stage.}\label{fig:expected}
\end{figure}

\begin{figure}[ht]
 \centering
     \setlength{\fboxrule }{0pt}
        \setlength{\fboxsep}{4pt}
\subfigure[Low arrival rate]{
        \fbox{\includegraphics[width=0.23\textwidth]{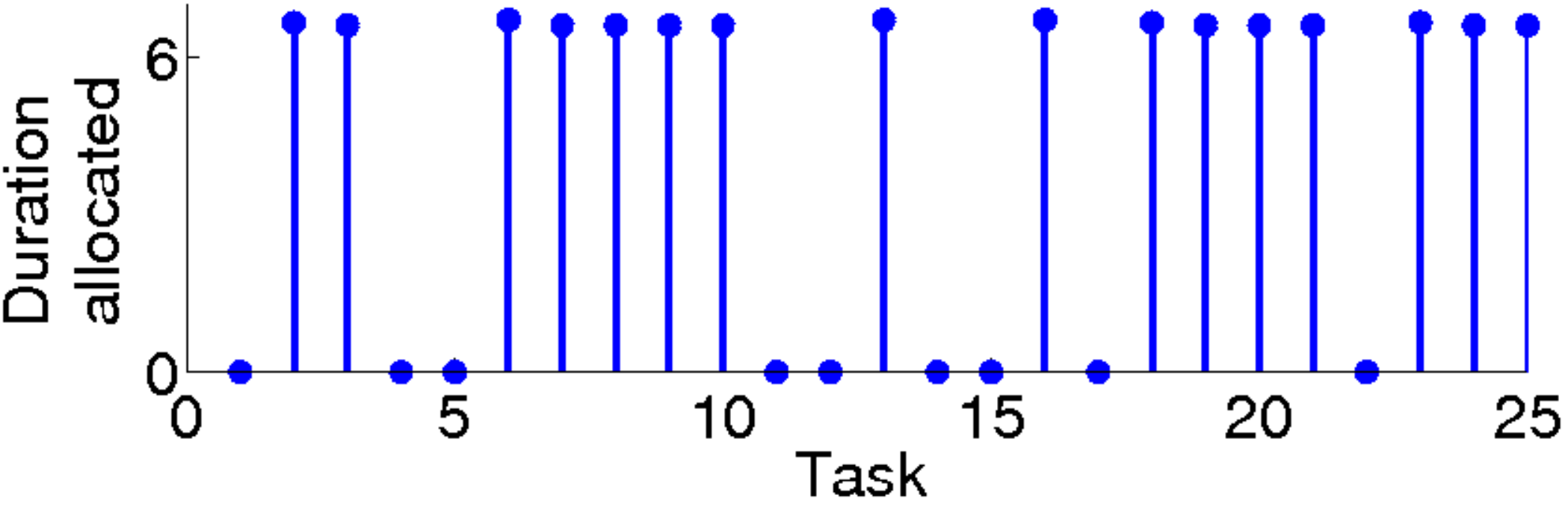}}
\fbox{\includegraphics[width=0.23\textwidth]{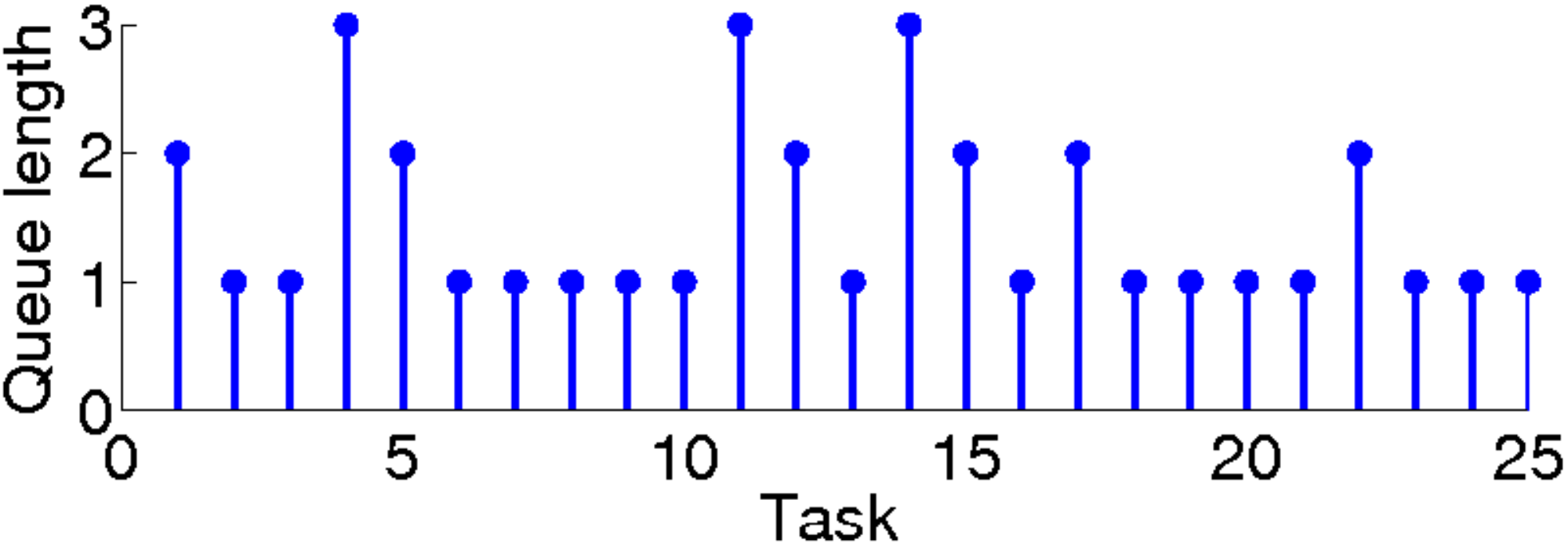}} }
\subfigure[Moderate arrival rate]{
        \fbox{\includegraphics[width=0.23\textwidth]{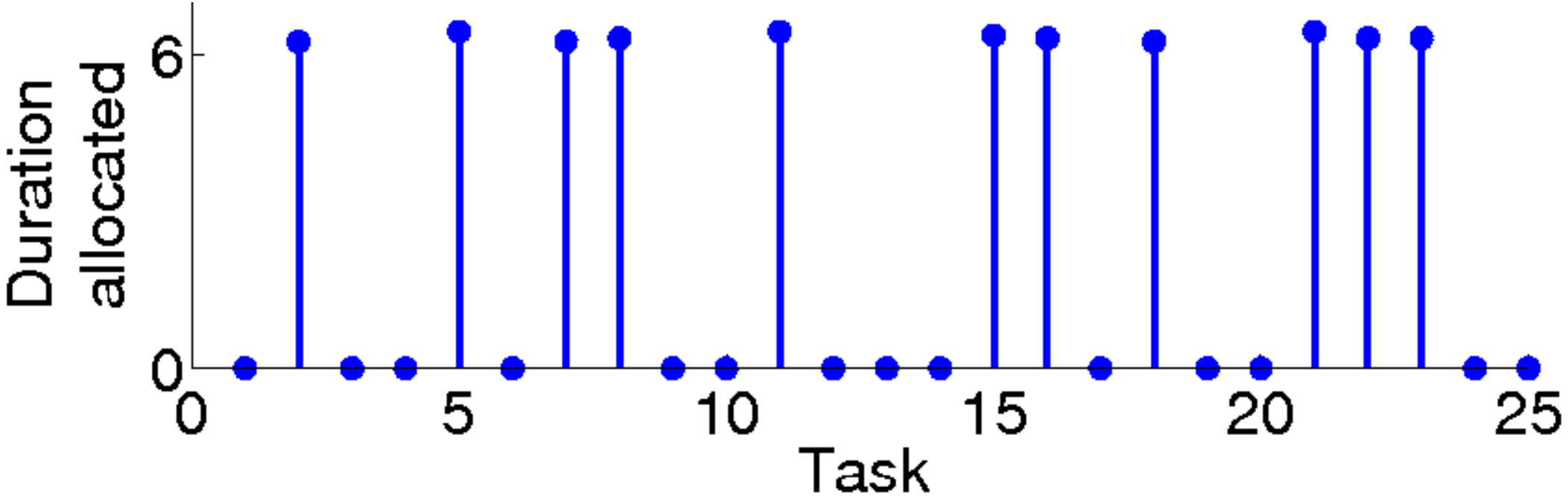}}
\fbox{\includegraphics[width=0.23\textwidth]{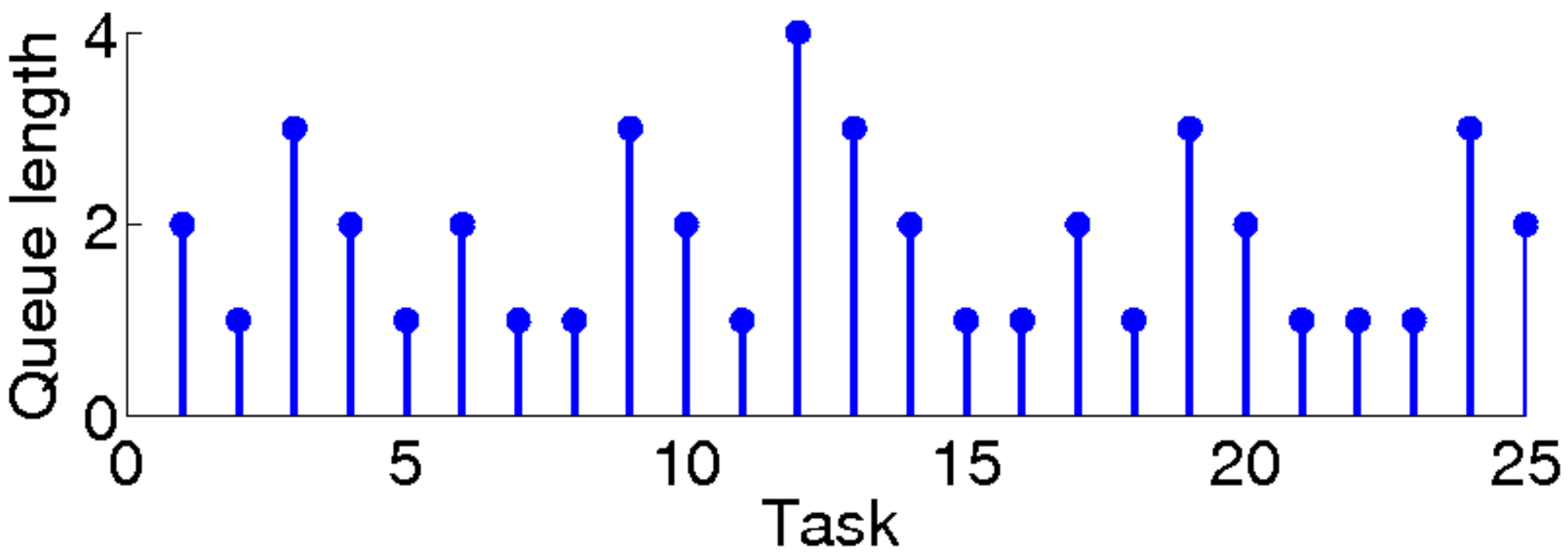}} }
\subfigure[High arrival rate]{
        \fbox{\includegraphics[width=0.23\textwidth]{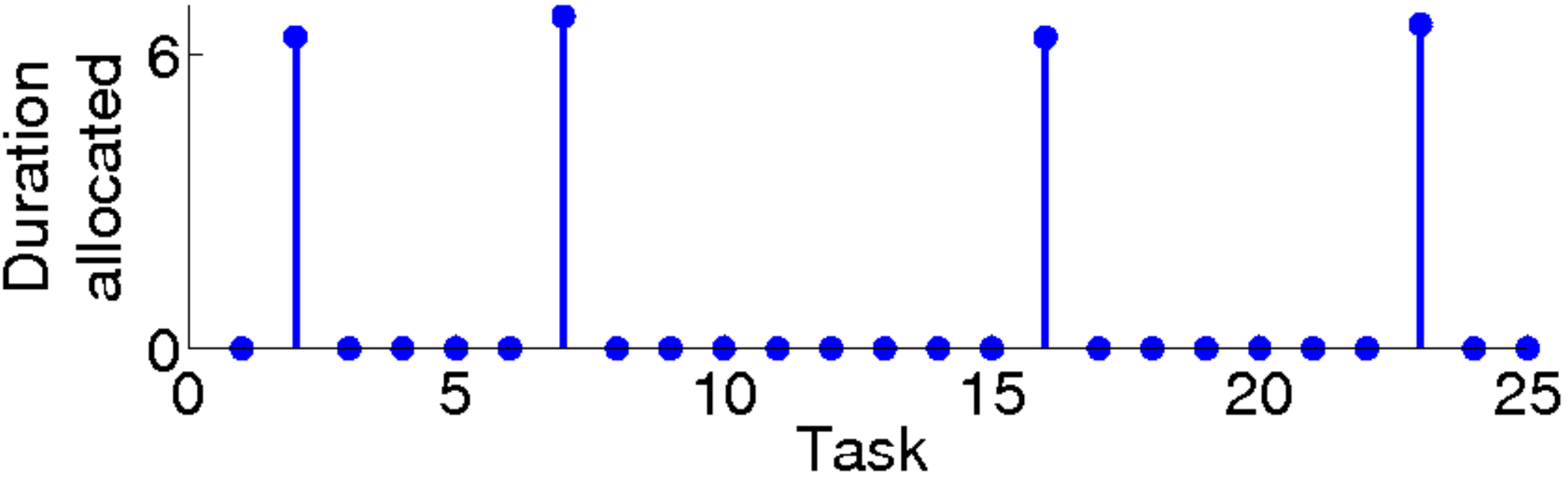}}
\fbox{\includegraphics[width=0.23\textwidth]{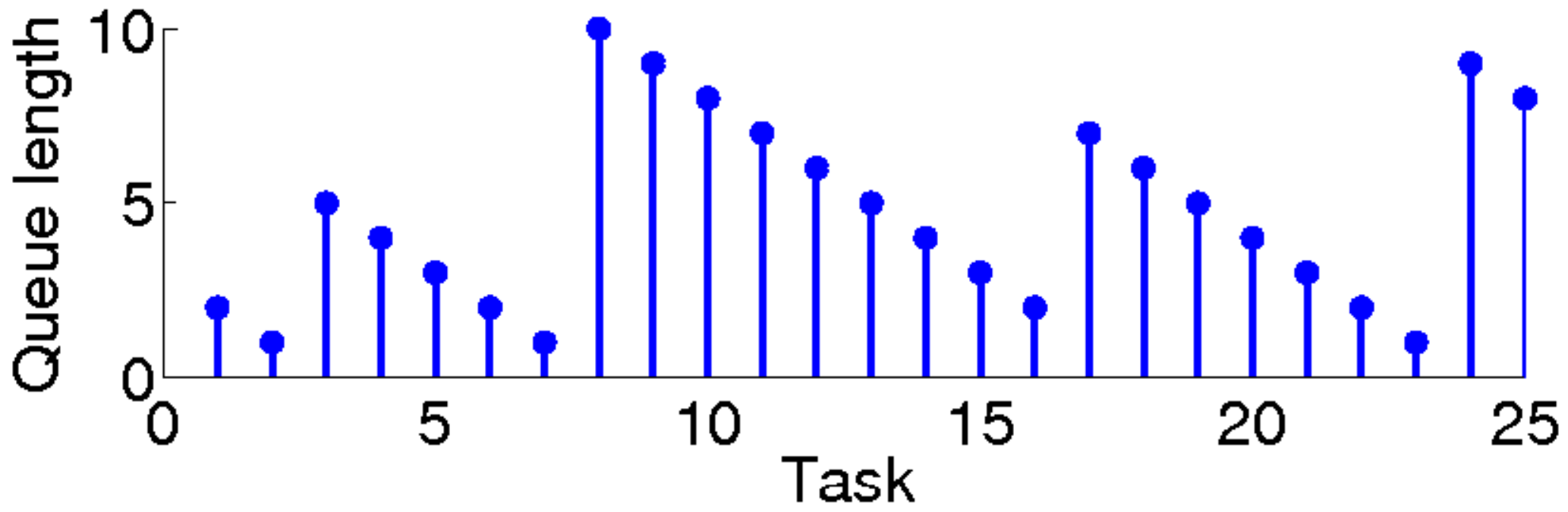}} }
\caption{Receding horizon policy for a {\it sample evolution} of the dynamic queue with latency penalty. An optimization problem with horizon length $N=10$ is solved at each stage.}\label{fig:sample}
\end{figure}

\begin{figure}[ht]
 \centering
     \setlength{\fboxrule }{0pt}
        \setlength{\fboxsep}{4pt}
        \fbox{\includegraphics[width=0.22\textwidth]{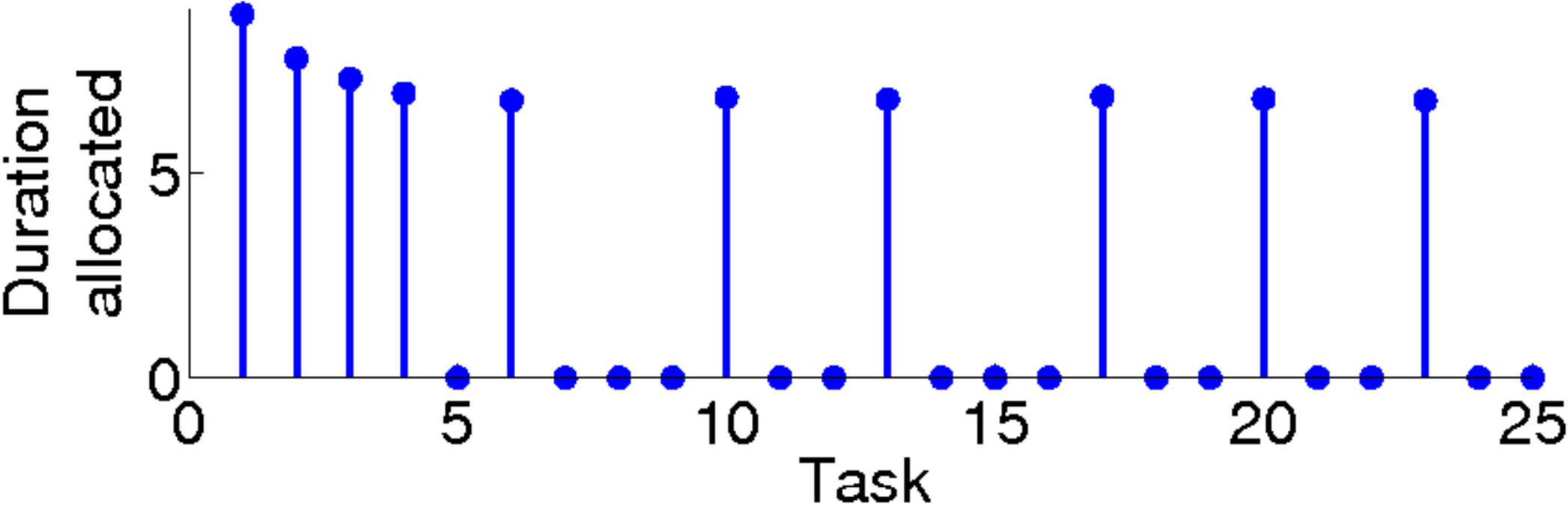}}
\fbox{\includegraphics[width=0.22\textwidth]{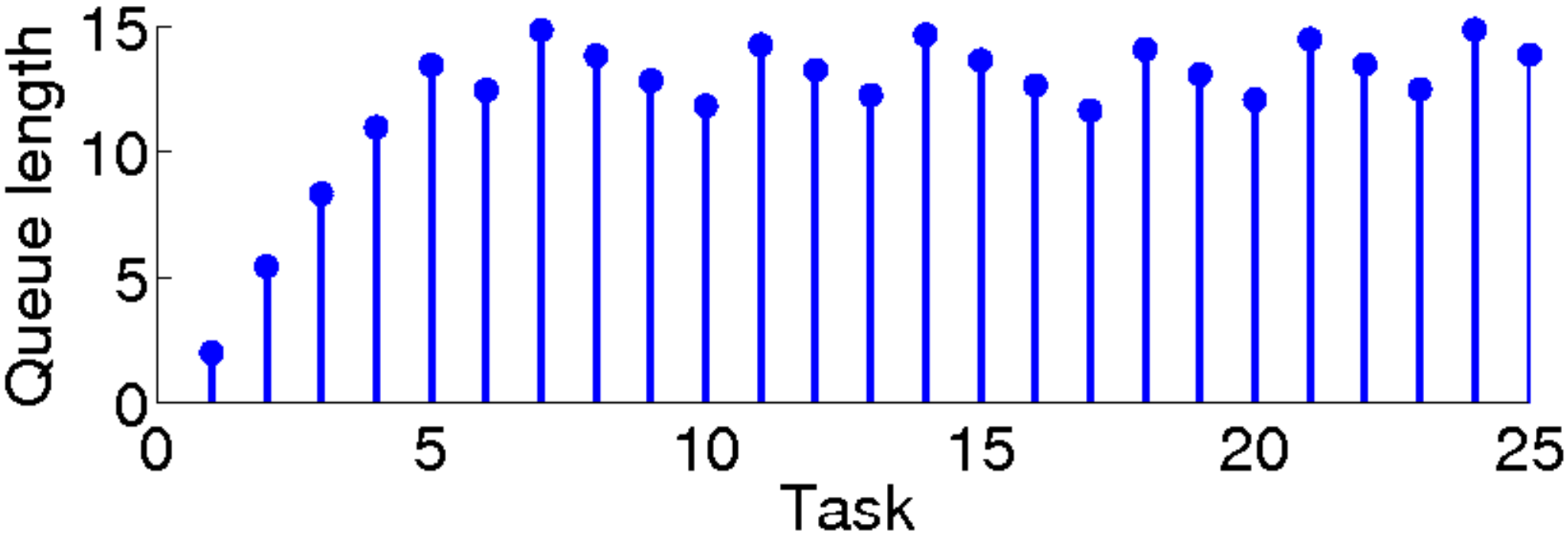}} 
\caption{Greedy policy for the expected evolution of the dynamic queue with latency penalty. An optimization problem with horizon length $N=1$ is solved at each stage.}\label{fig:greedy}
\end{figure}

\begin{figure}[ht]
 \centering
     \setlength{\fboxrule }{0pt}
        \setlength{\fboxsep}{4pt}
        \fbox{\includegraphics[width=0.22\textwidth]{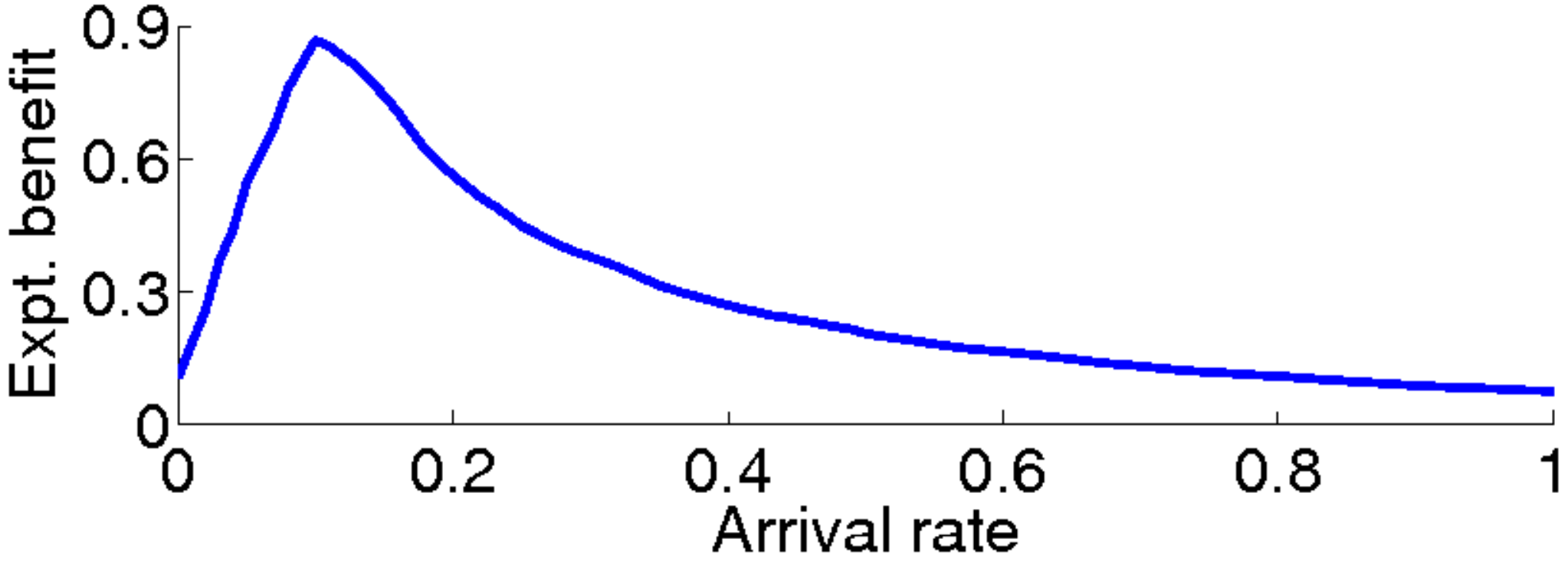}}
\fbox{\includegraphics[width=0.22\textwidth ]{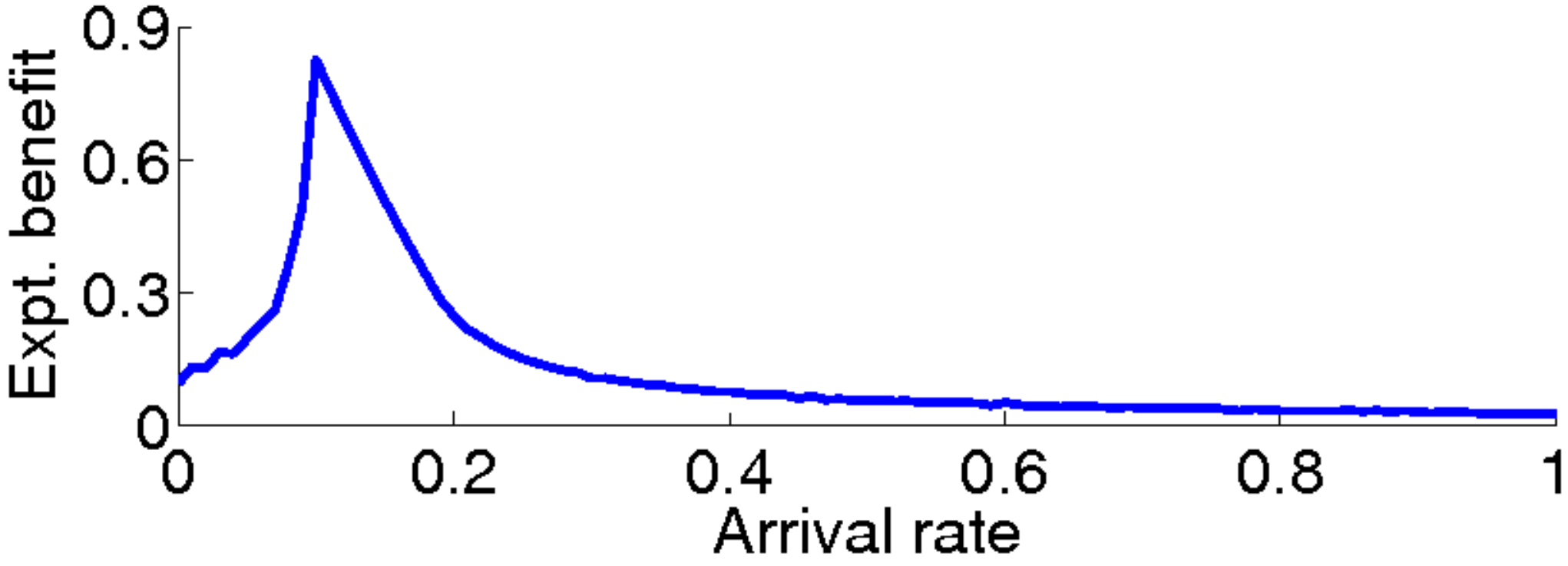}} 
\caption{Expected benefit over a finite horizon for the receding horizon and greedy policies}\label{fig:benefit}
\end{figure}

\begin{example}
If the human operator has to serve a queue of tasks with Poisson arrival at rate $\lambda=0.5$ per sec and the human receives
 an expected reward $f(t)=1/(1+\exp(5-t))$ for an allocation of duration $t$ secs to a task, while she incurs a penalty  $c=0.01$ 
 per sec for each task in queue. We solved an optimization problem with horizon length $N=10$ at each stage to determine the receding horizon optimization solution. 
  A receding horizon  policy for the expected evolution of the queue, at different arrival rates, is shown in Figure~\ref{fig:expected}.  
 A receding horizon duration allocation policy for a sample evolution of the queue, at different arrival rates, is shown in Figure~\ref{fig:sample}.
 The duration allocations for a greedy policy, i.e., an optimization with horizon length $N=1$ at each stage,  are shown in Figure~\ref{fig:greedy}.  The optimal expected benefit $J(\t)$ for the optimal and the greedy policy is shown in Figure~\ref{fig:benefit}. It can be seen that the maximum benefit is obtained at an arrival rate at which one expects only one task in the queue at each time. As expected, the performance of the greedy policy and the optimal policy is almost the same at this arrival rate. 
\oprocend
\end{example}

\begin{remark}[Optimal arrival rate]
It can be seen that the total expected benefit is maximum when there is always only one task in the queue. If there is more than one task in the queue, then the operator is incurring more penalty for the same reward. Thus, the optimal arrival is the one at which a task arrival is expected at the time when the previous task loses all its value, i.e., at $\tau^*=\max\setdef{t\in\real_{\ge 0}}{f'(t)=2c}$.  In general, there would be performance goals for the operator, and higher task arrival rate for the queue could be designed. Such a problem can be analyzed by putting a saturation on the sigmoidal performance function, and thus obtaining a new sigmoidal performance function.\oprocend
\end{remark}

\section{Conclusions }\label{sec:conclusion}
We presented an analysis on the decision making queues. Three particular problems were discussed. First, a time constrained decision making queue with no arrival was considered. We showed that the optimal policy may drop some tasks and assign equal time to the remaining tasks. Second, a decision making queue with no arrival and a latency penalty was considered. It was observed that the optimal policy may still drop some  tasks. Further, the duration allocation to the tasks increased with the decreasing queue length. Last, a decision making queue with latency penalty was considered. A receding horizon policy was developed to determine the optimal duration allocation. It was observed that the optimal policy may drop some tasks.

The decision support system designed in this paper assumes that the arrival rate of the tasks as well as the parameters in the performance function are known.  An interesting open problem is to estimate the parameters of the performance function, and come up with policies which perform an online estimation of the arrival rate, and accordingly determine the optimal allocation policy.

\end{document}